\documentclass{article}
\usepackage{amsfonts,amssymb,latexsym}

\newtheorem{theorem}{Theorem}[section]
\newtheorem{lemma}[theorem]{Lemma}
\newtheorem{corollary}[theorem]{Corollary}
\newtheorem{proposition}[theorem]{Proposition}
\newtheorem{conjecture}[theorem]{Conjecture}
\newenvironment{proof}
{\par\addvspace{0.3cm}\noindent{\rm Proof. }}
{\nopagebreak\mbox{}\hfill $\Box$\par\addvspace{0.25cm}}
\newenvironment{proofof}[1]
{\par\addvspace{0.3cm}\noindent{\rm Proof of #1. }}
{\nopagebreak\mbox{}\hfill $\Box$\par\addvspace{0.25cm}}

\newcommand{\C}{\mathbb{C}}
\newcommand{\Z}{\mathbb{Z}}

\newcommand{\R}{\mathbb{R}}
\newcommand{\T}{\mathbb{T}}

\newcommand{\cD}{\mathcal{D}}
\newcommand{\cL}{\mathcal{L}}
\newcommand{\cG}{\mathcal{G}}
\newcommand{\cO}{\mathcal{O}}
\newcommand{\cR}{\mathcal{R}}

\newcommand{\B}{B^1_1}
\newcommand{\G}{G_1}

\newcommand{\Ba}{\mathbf{a}}
\newcommand{\Bb}{\mathbf{b}}
\newcommand{\Bc}{\mathbf{c}}
\newcommand{\Bh}{\mathbf{h}}

\newcommand{\Bxi}{\mathbf{\xi}}
\newcommand{\Bet}{\mathbf{\eta}}
\newcommand{\Bnu}{\mathbf{\nu}}
\newcommand{\Bsi}{\mathbf{\sigma}}
\newcommand{\Bom}{\mathbf{\omega}}

\newcommand{\diag}{\mathrm{diag}\,}
\renewcommand{\Re}{\mathrm{Re}\,}

\newcommand{\sym}{\mathrm{sym}}

\newcommand{\subtwo}[2]%
{\ba{c}\scriptstyle #1\\\scriptstyle #2\ea}
\newcommand{\subthree}[3]%
{\ba{c}\scriptstyle #1\\\scriptstyle #2\\\scriptstyle #3\ea}

\newcommand{\bqn}{\begin{eqnarray}}
\newcommand{\eqn}{\end{eqnarray}}
\newcommand{\ba}{\begin{array}}
\newcommand{\ea}{\end{array}}
\newcommand{\nn}{\nonumber}

\newcommand{\si}{\sigma}
\newcommand{\ro}{\varrho}
\newcommand{\Ga}{\Gamma}

\newcommand{\Dp}{\cD_+'}
\newcommand{\Dm}{\cD_-'}
\newcommand{\Dpm}{\cD_\pm'}
\newcommand{\el}{\ell^2}
\newcommand{\iv}{^{-1}}
\newcommand{\iy}{\infty}

\newcommand{\wt}[1]{\widetilde{#1}}
\newcommand{\wh}[1]{\widehat{#1}}
\newcommand{\Cit}[1]{C^\iy(\T\setminus #1)}
\newcommand{\Cti}[1]{C^\iy_{#1}(\T)}
\newcommand{\Zm}{\Z_{-}}
\newcommand{\Zmn}{\Z_{-}\cup\{0\}}
\newcommand{\ros}{\ro_0,\ro_1,\ro_2}
\newcommand{\mi}{{\mu_1}}
\newcommand{\mii}{{\mu_2}}
\newcommand{\refeq}[1]{{\rm (\ref{#1})}}

\begin{document}

\date{}
\title{Asymptotic formulas for the determinants of
symmetric Toeplitz + Hankel matrices}
\author{Estelle L. Basor\thanks{ebasor@calpoly.edu.
          Supported in part by NSF Grant DMS-9970879.}\\
               Department of Mathematics\\
               California Polytechnic State University\\
               San Luis Obispo, CA 93407, USA
        \and
        Torsten Ehrhardt\thanks{tehrhard@mathematik.tu-chemnitz.de.}\\
                Fakult\"{a}t f\"{u}r Mathematik\\
                Technische Universit\"{a}t Chemnitz\\
                09107 Chemnitz, Germany}
\maketitle

\begin{abstract}
We establish asymptotic formulas for the determinants
of $N\times N$ Toeplitz + Hankel matrices $T_N(\phi)+H_N(\phi)$ as $N$ goes
to infinity for singular generating functions $\phi$ defined on the unit
circle in the special case where $\phi$ is even, i.e.,
where the Toeplitz + Hankel matrices are symmetric.
\end{abstract}


\section{Introduction}

In the theory of random matrices, for certain ensembles,
one is led to consider the
asymptotics of Fredholm operators of the form $I+W+H$ where $W$ is a
finite and symmetric Wiener-Hopf operator and $H$ is a finite Hankel
operator \cite{Me}.  This problem arises when investigating the probability
distribution function of a random variable thought of as a function of
the eigenvalues of a positive Hermitian random matrix.  For general
information about random matrix theory we refer the reader to
\cite{Me} and also to
\cite{F,B,BT} for more specific tie-ins to the random variable
problem.

The focus of this paper is to study the discrete analogue of this
problem.  This is not precisely the desired situation for those
interested in random matrix theory.  However, it is a natural starting
place for cases where the random variable is discontinuous, since then
the discrete nature of the computations make things a bit more
accessible and the mathematical questions that arise are quite
interesting in themselves.

The discrete analogue of this problem is to find an asymptotic expansion
of the determinants of Toeplitz + Hankel matrices
\bqn\label{f1.1}
M_N(\phi) &=& T_N(\phi)+H_N(\phi)
\eqn
in the case where these matrices are symmetric.
Here the $N\times N$ Toeplitz and Hankel matrices are
defined as usual by
\begin{equation}\label{f1.2}
T_N(\phi) = \big(\phi_{j-k}\big)_{j,k=0}^{N-1},\qquad\quad
H_N(\phi) = \big(\phi_{j+k+1}\big)_{j,k=0}^{N-1}.
\end{equation}
The entries $\phi_n$ are the Fourier coefficients
\bqn
\phi_n &=&
\frac{1}{2\pi}\int_0^{2\pi}\phi(e^{i\theta})e^{-in\theta}\,d\theta
\eqn
of a function $\phi\in L^1(\T)$ defined on the unit circle $\T$.

The matrices $M_N(\phi)$ are symmetric if and only if the function
$\phi$ is even, i.e., if $\phi(e^{i\theta})=\phi(e^{-i\theta})$.
{}From the point of view of random matrix theory, one is particularly
interested in the asymptotics of $\det M_N(\phi)$ as $N\to\iy$ in the
case of even, piecewise continuous functions $\phi$. This is related
to the problem of finding the  distribution function of a random
variable that counts the number of eigenvalues of a random matrix that
lie in an interval and to finding the distribution function for other
random variables. See \cite{B} to see the connections between these
problems.

The problem of determining the asymptotics of the determinants of (not
necessarily symmetric) matrices $M_N(\phi)$ has been studied intensively in
a previous paper \cite{BE1}.  For example, it was shown there
that if $\phi$ is continuous and sufficiently smooth, then the
asymptotics are very similar to the ones given in the Strong Szeg\"{o}-Widom
Limit Theorem.  Indeed, it is only in the constant, or third
order term that the answers differ.  This is no surprise since if
$\phi$ is continuous, then the Toeplitz operator is perturbed by a
compact Hankel operator only.

If, however, the symbol $\phi$ is singular, then the problem is much
harder to solve.  In the case of Toeplitz determinants the answer is
provided by the Fisher-Hartwig conjecture, which has been proved
under certain smoothness assumptions in all the cases where it is
expected to hold \cite{E1}.  In \cite{BE1} an asymptotic formula for
the determinants $\det M_N(\phi)$ was obtained for piecewise
continuous functions $\phi$, but under the additional assumption that
the function $\phi$ does not possess a discontinuity at both a point
on the unit circle and its complex conjugate.  In this case, the
asymptotic formula shows that the asymptotics differ from the
asymptotics of Toeplitz determinants only in the third order term,
i.e., in the constant term, while the second order term is the same.

Unfortunately, the additional assumption on the location of the
discontinuities imposed in \cite{BE1} excludes all even, piecewise
continuous functions.  Hence the paper \cite{BE1} does not
answer the discrete analogue of the problem motivated by random
matrix theory.

It is the purpose of this paper to solve this problem by establishing
an asymptotic formula for determinants of matrices (\ref{f1.1}) for
even piecewise continuous functions $\phi$.  For such functions, it
turns out that the asymptotics differs also in the second order term
in comparison with the asymptotics of Toeplitz determinants.

The paper is organized as follows. In Section \ref{sec:2} we will recall
some of the results established in \cite{BE1} that are of relevance for this
paper. In Section \ref{sec:3} we establish an identity which is the key for
computing the asymptotics of $\det M_N(\phi)$ for $\phi$ even.
This identity can be formulated as follows:
\bqn\label{f1.4new}
(\det M_N(\Ba))^2&=& \det T_{2N}(\Ba\Bsi).
\eqn
In this identity $\Ba$ can no longer be considered as a function, but has
to be understood as a distribution, which satisfies certain properties.
Moreover,
$\Bsi$ is here a certain concrete distribution. This identity was
established
first in \cite{Kra} in a formulation that is not based on distributions.
The goal of Section \ref{sec:3} is to provide the necessary tools needed for
dealing with distributions, in particular, to define a product between $\Ba$
and
$\Bsi$ in an appropriate way. Having done this, we are able to derive the
distributional formulation of this identity from the original one.

The identity (\ref{f1.4new}) reduces the asymptotics of $\det M_N(\Ba)$
to the asymptotics of the (skewsymmetric) Toeplitz determinant
$\det T_{2N}(\Ba\Bsi)$. In order to analyze this Toeplitz determinant we
cannot
rely on the (original) Fisher-Hartwig conjecture because it breaks down in
this
case. However, in Section \ref{sec:4} we will proof a limit theorem saying
that
the quotient
\bqn
\frac{\det T_{2N}(\Ba\Bsi)}{\det T_{2N}(\Ba)}
\eqn
converges -- under certain conditions on $\Ba$ -- to a nonzero constant.
In order to prove this limit theorem we make heavy use of the machinery that
has been developed in \cite{E1} in order to prove the Fisher-Hartwig
conjecture.

Thus, up to this point, we have reduced the asymptotics of $\det M_N(\Ba)$
to the asymptotics of $\det T_{2N}(\Ba)$. The Toeplitz determinant
$\det T_{2N}(\Ba)$ is (generically) of a kind for which the Fisher-Hartwig
conjecture holds. In Section \ref{sec:5} we will therefore recall the
Fisher-Hartwig conjecture in the form as it has been proved in \cite{E1}.
Moreover, we specialize it to the distributions (namely, even distributions
of
Fisher-Hartwig type) that we are interested in. In Section
\ref{sec:6} we combine all the previous results and obtain the
asymptotics of $\det M_{N}(\Ba)$ for even distributions $\Ba$ of
Fisher-Hartwig
type, which satisfy appropriate conditions on the parameters.

In Section \ref{sec:7} we finally specialize the quite general result of
Section \ref{sec:6} to even piecewise continuous functions. We thus obtain
the
asymptotics of $\det M_{N}(\phi)$ for a certain class of even piecewise
continuous
functions. This result together with the results that are known from
\cite{BE1} suggest a conjecture about the asymptotics
of $\det M_{N}(\phi)$ for quite general piecewise continuous, not just
those necessarily even. We end the paper with the conjecture.


\section{Known results for determinants of Toeplitz plus Hankel matrices}
\label{sec:2}

Let us begin by recalling some of the results concerning the asymptotics of
the
determinants of the matrices $M_N(\phi)$ that have already been
established in \cite{BE1}.

We first consider the case of continuous and sufficiently smooth
generating functions $\phi$, where an analogue to the Strong Szeg\"o-Widom
Limit Theorem holds.  In order to be more specific about the
smoothness condition, let us consider the Besov class $\B$, which is
by definition the set of all functions $b\in L^1(\T)$ such that
\bqn
\|b\|_{\B} &:=&
\int_{-\pi}^\pi\frac{1}{y^2}\int_{-\pi}^\pi
\left|b(e^{ix+iy})+b(e^{ix-iy})-2b(e^{ix})\right|\,dxdy<\iy.
\eqn
It is known that $\B$ forms a Banach algebra with the above norm and
is continuously embedded into the Banach algebra of all
continuous functions on $\T$.  By $\G\B$ we denote
the set of all nonvanishing functions in $\B$ with winding number
zero.  The set $\G\B$ can also be characterized as the set of all
functions $b$ which possess a logarithm $\log b$ in $\B$.

For $b\in \G\B$, the constants
\bqn
G[b] &=& \exp\left([\log b]_{0}\right),
 \label{f1.Gb}\\
E[b] &=&\exp\left(\sum_{k=1}^\iy k[\log b]_k[\log b]_{-k}\right),
 \label{f1.Eb}\\
F[b] &=&
 \exp\left(\sum_{k=1}^\iy [\log b]_{2k-1}-
 \frac{1}{2}\sum_{k=1}^\iy k[\log b]_k^2\right),
 \label{f1.Fb}
\eqn
are well defined, where $[\log b]_{n}$ stand for the Fourier
coefficients of $\log b\in \B$.
Moreover, for $b\in\G\B$, the functions $b_+,b_-\in\G\B$ are
well-defined by
\bqn\label{f1.bpm}
b_\pm(t) &=& \exp\left(\sum_{n=1}^\iy t^{\pm n}[\log b]_{\pm n}\right),
\qquad t\in\T.
\eqn
Note that $b(t)=b_-(t)G[b]b_+(t)$, $t\in\T$, is the normalized canonical
Wiener-Hopf factorization of the function $b$.

The analogue to the Strong Szeg\"o-Widom Limit Theorem for the determinants
$\det M_N(\phi)$, which has been established in \cite[Corollary 2.6]{BE1},
now says that if $b\in \G\B$, then
\bqn
\det M_N(b) &\sim&
G[b]^NE[b]F[b]\quad\mbox{ as }N\to\iy.
\eqn
In the case of even functions $b\in\G\B$ this simplifies to
\bqn
\det M_N(b) &\sim&
G[b]^N\wh{E}[b]\quad\mbox{ as }N\to\iy,\label{f2.7new}
\eqn
where $\wh{E}[b]$ is the constant
\bqn\label{f1.Ehb}
\wh{E}[b] &=& \exp\left(\frac{1}{2}\sum_{k=1}^\iy k[\log b]_k^2
+\sum_{k=1}^\iy [\log b]_{2k-1}\right).
\eqn

In order to discuss the asymptotics for the case of piecewise continuous
generating
functions $\phi$, let us introduce the functions
\bqn
t_{\beta,\theta_0}(e^{i\theta}) &=&
e^{i\beta(\theta-\theta_0-\pi)},\qquad 0<\theta-\theta_0<2\pi,
\eqn
where $\beta\in\C$ and $\theta_0\in(-\pi,\pi]$. The piecewise
continuous functions that we consider are of the form
\bqn\label{f1.9z}
\phi(e^{i\theta}) &=& b(e^{i\theta}) \prod_{r=1}^R
t_{\beta_r,\theta_r}(e^{i\theta}),
\eqn
where $\theta_1,\dots,\theta_R\in(-\pi,\pi]$ are distinct numbers
determining the location of the jump discontinuities and
$\beta_1,\dots,\beta_R$ are complex parameters determining the
``size'' of the jumps.  The function $b$ is usually assumed to belong
to $\G\B$.

As is known from the theory of Toeplitz determinants, a key ingredient
for the determination of the asymptotics of the determinants are
localization theorems (see, e.g., \cite{BS}).
A localization theorem for the determinants
$\det M_N(\phi)$ with $\phi$ being piecewise continuous has been
established in \cite[Theorem~5.11]{BE1}.  This localization theorem
reduces the asymptotics of $\det M_N(\phi)$ for ``general'' piecewise
continuous functions (\ref{f1.9z}) to the asymptotics for particular
piecewise continuous functions.

\begin{theorem}[Localization Theorem]\label{t2.1new}
Let $\phi$ be a function of the form
\bqn
\phi(e^{i\theta}) &=&
b(e^{i\theta}) \phi_+ (e^{i\theta}) \phi_-(e^{i\theta})
\prod_{r=1}^R \phi_r(e^{i\theta}),\label{f4.LocForm}
\eqn
where $b\in \G\B$, $\phi_+=t_{\beta_+,0}$, $\phi_-=t_{\beta_-,\pi}$
and $\phi_r=t_{\beta^+_r,\theta_r}t_{\beta^-_r,-\theta_r}$
for $1\le r\le R$.
Suppose that $\theta_1,\dots,\theta_R\in(0,\pi)$ are distinct numbers and
that $\beta_\pm,\beta_1^\pm,\dots,\beta_R^\pm\in\C$ are such that
\begin{itemize}
\item[(a)]
$-1/2<\Re\beta_+<1/4$ and $-1/4<\Re\beta_-<1/2$,
\item[(b)]
$|\Re\beta_r^+|<1/2$ and $|\Re\beta_r^-|<1/2$ and
$|\Re(\beta_r^++\beta_r^-)|<1/2$ for each $1\le r\le R$.
\end{itemize}
Then
\bqn
\lim_{N\to\iy} \frac{\det M_N(\phi)}{\displaystyle
\det M_N(b) \det M_N(\phi_+ ) \det M_N(\phi_-) \prod_{r=1}^R
\det M_N(\phi_r)} &=& H,\nn
\eqn
where
\bqn
H &=&
b_+(1)^{2\beta_+}b_-(1)^{-\beta_+}
b_+(-1)^{2\beta_-}b_-(-1)^{-\beta_-}
2^{3\beta_+\beta_-}
\nn\\
&&\times\prod_{r=1}^R
b_+(t_r)^{\beta^+_r+\beta^-_r}
b_-(t_r)^{-\beta^+_r}
b_+(t_r\iv)^{\beta^+_r+\beta^-_r}
b_-(t_r\iv)^{-\beta^-_r}
\nn\\
&&\times\prod_{r=1}^R
(1-t_r)^{\beta_+(\beta^+_r+2\beta^-_r)}
(1-t_r\iv)^{\beta_+(2\beta^+_r+\beta^-_r)}
\nn\\
&&\times\prod_{r=1}^R
(1+t_r)^{\beta_-(\beta^+_r+2\beta^-_r)}
(1+t_r\iv)^{\beta_-(2\beta^+_r+\beta^-_r)}
\nn\\
&&\times\!\!\prod_{1\le r<s\le R}\!\!
(1-t_rt_s)^{\beta^-_r\beta^-_s+\beta^+_r\beta^-_s+\beta^-_r\beta^+_s}
(1-t_r\iv t_s\iv)^{\beta^+_r\beta^+_s+\beta^+_r\beta^-_s+\beta^-_r\beta^+_s}
\nn\\
&&\times\!\!\prod_{1\le r<s\le R}\!\!
(1-t_rt_s\iv)^{\beta^+_r\beta^+_s+\beta^-_r\beta^-_s+\beta^-_r\beta^+_s}
(1-t_r\iv t_s)^{\beta^+_r\beta^+_s+\beta^-_r\beta^-_s+\beta^+_r\beta^-_s}.
\nn
\eqn
Here $t_r=e^{i\theta_r}$, $1\le r\le R$, and $b_\pm$ are the functions
\refeq{f1.bpm}.
\end{theorem}

The asymptotic behavior of $\det M_N(\phi)$ with the generating
function $\phi_+=t_{\beta_+,0}$ and $\phi_-=t_{\beta_-,\pi}$,
respectively, has also been determined (see Theorem 6.2 and Theorem
6.3 in \cite{BE1}).  These functions have one single jump discontinuity
at the points $1$ and $-1$, respectively.

\begin{theorem}\label{t2.2new}
Let $\beta\in\C\setminus\Z$. Then
\begin{itemize}
\item[(a)]
$\displaystyle\lim_{N\to\iy}
\frac{\det M_N(t_{\beta,0})}{N^{-\frac{3\beta^2}{2}-\frac{\beta}{2}}} =
(2\pi)^{\frac{\beta}{2}}2^{\frac{3\beta^2}{2}}
\frac{G(\frac{1}{2}-\beta)G(1-\beta)G(1+\beta)}{G(\frac{1}{2})}$,
\item[(b)]
$\displaystyle\lim_{N\to\iy}
\frac{\det M_N(t_{\beta,\pi})}{N^{-\frac{3\beta^2}{2}+\frac{\beta}{2}}} =
(2\pi)^{\frac{\beta}{2}}2^{\frac{3\beta^2}{2}}
\frac{G(\frac{3}{2}-\beta)G(1-\beta)G(1+\beta)}{G(\frac{3}{2})}$.
\end{itemize}
\end{theorem}

In these asymptotic formulas the Barnes $G$-function $G(z)$
appears \cite{Bar,WW}, which is an entire function defined by
\bqn
 G(1+z) &=&(2\pi)^{\frac{z}{2}}e^{-\frac{(z+1)z}{2}-C_{E}\frac{z^{2}}{2}}
 \prod_{k=1}^\iy\left(\left(1+\frac{z}{k}\right)^{k}
 e^{-z+\frac{z^{2}}{2k}}\right)\label{9}
\eqn
with $C_{E}$ being Euler's constant.

In the case of the generating functions
$\phi_r=t_{\beta^+_r,\theta_r}t_{\beta^-_r,-\theta_r}$, which have two
jump discontinuities at a point of the unit circle and its complex
conjugate, the asymptotic behavior is only known in particular cases.
One case is that where either $\beta^+_r=0$ or $\beta^-_r=0$, i.e.,
where the function has exactly one jump discontinuity at a point in
$\T\setminus\{1,-1\}$.  Here the result is taken from
\cite[Theorem 4.5]{BE1}.

\begin{theorem}\label{t2.3new}
Let $\theta_0\in(-\pi,0)\cup(0,\pi)$ and $\beta\in\C$ be such that
$|\Re\beta|<1/2$. Put $t_0=e^{i\theta_0}$. Then
\bqn
\lim_{N\to\iy}\frac{\det M_N(t_{\beta,\theta_0})}{ N^{-\beta^2}} &=&
G(1-\beta)G(1+\beta)
\left(1-t_0\iv\right)^{\frac{\beta^2}{2}+\frac{\beta}{2}}
\left(1+t_0\iv\right)^{\frac{\beta^2}{2}-\frac{\beta}{2}}.\nn
\eqn
\end{theorem}

Finally, another, even more particular case of a function with two jump
discontinuities at $i$ and $-i$ and the same size of the jumps
has been established \cite{BE1} if one combines Theorem 7.4 and Theorem 7.5.

\begin{theorem}\label{t2.4new}
Let $\beta\in\C\setminus\Z$. Then
\bqn
\lim_{N\to\iy}
\frac{\det
M_N(t_{\beta,\frac{\pi}{2}}t_{\beta,-\frac{\pi}{2}})}{N^{-3\beta^2}}
&=& 2^{4\beta^2}G(1-2\beta)G(1+\beta)^2.\nn
\eqn
\end{theorem}

For later use, let us specialize the localization theorem
(Theorem \ref{t2.1new}) to the case of even functions $\phi$ which are of
the form
(\ref{f4.LocForm}).

\begin{corollary}\label{c2.5new}
Let $\phi$ be a function of the form
\bqn
\phi(e^{i\theta}) &=&
b(e^{i\theta})
\prod_{r=1}^R \phi_r(e^{i\theta}),
\eqn
where $b\in \G\B$ is an even function and
$\phi_r=t_{\beta_r,\theta_r}t_{-\beta_r,-\theta_r}$ for $1\le r\le R$.
Suppose that
$\theta_1,\dots,\theta_R\in(0,\pi)$ are distinct numbers and that
$\beta_1,\dots,\beta_R\in\C$ are such that $|\Re\beta_r|<1/2$
for each $1\le r\le R$. Then
\bqn
\lim_{N\to\iy} \frac{\det M_N(\phi)}{\displaystyle
\det M_N(b) \prod_{r=1}^R \det M_N(\phi_r)} &=& H,\nn
\eqn
where
\bqn
H &=&
\prod_{r=1}^R
b_+(t_r)^{\beta_r}
b_-(t_r)^{-\beta_r}
\nn\\
&&
\times\!\!\prod_{1\le r<s\le R}\!\!
(1-t_rt_s)^{-\beta_r\beta_s}
(1-t_r\iv t_s\iv)^{-\beta_r\beta_s}
\nn\\
&&\times\!\!\prod_{1\le r<s\le R}\!\!
(1-t_rt_s\iv)^{\beta_r\beta_s}
(1-t_r\iv t_s)^{\beta_r\beta_s}.\nn
\eqn
Here $t_r=e^{i\theta_r}$, $1\le r\le R$, and $b_\pm$ are the functions
\refeq{f1.bpm}.
\end{corollary}
\begin{proof}
We apply Theorem \ref{t2.1new} with the parameters
$\beta_\pm=0$ and $\beta_r^\pm=\pm\beta_r$. We remark also that
$b_+=\tilde{b}_-$.
\end{proof}


\section{Preliminary results for determinants of \\ symmetric Toeplitz plus
Hankel matrices}
\label{sec:3}

The first step in order to determine the asymptotics of the determinants
of symmetric Toeplitz + Hankel matrices $M_N(\phi)$ is to express
these determinants by means of determinants of skewsymmetric Toeplitz
matrices.
The identity as it appears in the following theorem
has been stated explicitly in \cite[Lemma~18]{Kra}, but is already
implicitly
contained in \cite[Lemma~1]{Gor} and
\cite[Proof of Theorem~7.1(a)]{Stem}, where it has been proved.
A different, self-contained proof has been given by the authors in
\cite{BE2}.

\begin{theorem}\label{t2.1}
Let $\{a_n\}_{n=-\iy}^{\iy}$ be a sequence of complex numbers such that \\
$a_{-n}=a_{n}$. Let
\bqn\label{f.cn}
c_{n} &=& \sum_{k=-n+1}^n a_{k}\qquad \mbox{ for }n>0,
\eqn
and put $c_0=0$ and $c_{-n}=-c_{n}$. Then
\bqn\label{f2.2}
\left(\det \left(a_{j-k}+a_{j+k+1}\right)_{j,k=0}^{N-1}\right)^2
&=&
\det\left(c_{j-k}\right)_{j,k=0}^{2N-1}.
\eqn
\end{theorem}

The matrices appearing in (\ref{f2.2}) are a symmetric Toeplitz + Hankel
matrix
of the same kind as (\ref{f1.1})
and a skewsymmetric Toeplitz matrix. If we are trying to rewrite this
identity by using the standard notation (\ref{f1.2}) for Toeplitz
and Hankel matrices where the sequences $\{a_n\}$ and $\{c_n\}$
are the Fourier coefficients of functions $a,c\in L^1(\T)$, we face
the difficulty that this is in general not possible. Consider for instance
the
simplest case where $a(t)=1$, i.e., $a_0=1$ and $a_n=0$ if $n\neq0$.
Then we obtain $c_n={\rm sign}(n)$, and obviously, there does not exist
a function $c\in L^1(\T)$ with such Fourier coefficients.

A way out of this situation is to consider distributions on the unit circle
in place of functions in $L^1(\T)$ and to take their Fourier coefficients
as the entries of the Toeplitz and Hankel matrices.
For this purpose we need several preliminary results. Apart from
basic issues, the following has all been stated in \cite{E2} and proved in
\cite{E1}.

Let $\cD=C^\iy(\T)$ be the linear topological space of all infinitely
differentiable
functions defined on the unit circle. By $\cD'$ we denote the set of all
distributions on the unit circle, i.e., linear and continuous functionals
on $\cD$. The Fourier coefficients of a distribution $\Ba\in\cD$
are defined as
\bqn
\Ba_n &=& \Ba(\chi_{-n}),
\eqn
where $\chi_n\in\cD$ is the function $\chi_n(t)=t^n$.
There is a natural identification of functions $a\in L^1(\T)$ with a subset
of
distributions. It is established by the mapping $a\in L^1(\T)\mapsto
\Ba\in\cD'$ where
\bqn\label{f2.4}
\Ba(f) &=& \frac{1}{2\pi}\int_{0}^{2\pi}a(e^{i\theta})f(e^{i\theta})
\,d\theta,\qquad f\in\cD.
\eqn
This definition ensures that the Fourier coefficients of $a$ and $\Ba$ are
the same. We also remark that there is a one-to-one correspondence between
$\cD'$ and the set of all at most polynomially increasing sequences
$\{a_n\}_{n=-\iy}^\iy$, which is given by associating
to $\Ba\in\cD'$ the series $\{\Ba_n\}_{n=-\iy}^\iy$
of its Fourier coefficients.

Let $a\in\cD$ and $\Bb\in\cD'$. Then the product of $a$ and $\Bb$
is the distribution $a\Bb\in\cD'$ which is defined by
\bqn
(a\Bb)(f) &=& \Bb(af), \qquad
f\in\cD.
\eqn

Let $K$ be a compact subset of $\T$. We denote by $\Cit{K}$ the set of all
infinitely differentiable functions on $\T\setminus K$. By $\Cti{K}$ we
refer to the set of all functions $f\in\cD$ which vanish on an open
neighborhood of $K$. The product of a function $f\in\Cit{K}$ with a function
$g\in\Cti{K}$ is a function $fg\in\Cti{K}\subseteq\cD$ by putting
$(fg)(t)=0$ for $t\in K$.

We will proceed with some definitions that are not so quite common,
but necessary for our considerations. They are taken from \cite{E1,E2}.
Let $\cD'(K)$ stand for the set of all distributions $\Ba$ for which there
exists a function $a\in\Cit{K}$ such that
\bqn\label{f1.9}
\Ba(f) &=& \frac{1}{2\pi}\int_0^{2\pi}a(e^{i\theta})f(e^{i\theta})\,d\theta
\eqn
for all $f\in\Cti{K}$. The function $a\in\Cit{K}$ is uniquely determined by
$\Ba$ and called the {\em smooth part of the distribution $\Ba$}.
Definition (\ref{f1.9}) can be rephrased by saying that
$f\Ba=fa$ for all $f\in\Cti{K}$, where the left hand side of this equation
is a distribution $\cD'$ and the right hand side is a function in $\cD$,
which are identified in the sense of (\ref{f2.4}).

Next we are going to show that one can define -- under certain assumptions
-- the product of two distributions. Let $M$ and $N$ be compact and disjoint
subsets of the unit circle. Given $\Ba\in\cD'(M)$ and $\Bb\in\cD'(N)$
with smooth parts $a\in\Cit{M}$ and $b\in\Cit{N}$, let
$\Ba\Bb\in\cD'(M\cup N)$ be defined as
\bqn\label{f1.10}
\Ba\Bb &=& (bf_b)\Ba+(af_a)\Bb,
\eqn
where $f_a\in\Cti{M}$ and $f_b\in\Cti{N}$ are such that $f_a+f_b=1$.
This definition is independent of the particular choice of $f_a$ and $f_b$.
Moreover, $\Ba\Bb$ has the smooth part $ab\in\Cit{(M\cup N)}$.

Given a function $a$ defined on (a subset of) the unit circle, we define
the function $\tilde{a}$ by $\tilde{a}(t)=a(t\iv)$, $t\in\T$.
In accordance with this definition, given a distribution $\Ba\in\cD'$,
let $\tilde{\Ba}\in\cD'$ stand for the distribution with Fourier
coefficients $\tilde{\Ba}_{-n}=\Ba_n$. A distribution $\Ba$ will be called
{\em even} ({\em odd}) if $\Ba=\pm \tilde{\Ba}$. A function $a$ will be
called
{\em even} ({\em odd}) if $a=\pm \tilde{a}$. Finally, if $K$ is a subset of
$\T$, put $\wt{K}=\{t\in\T:t\iv\in K\}$. If $\wt{K}=K$, we call $K$ a
{\em symmetric subset of the unit circle}.

In the reformulation of Theorem \ref{t2.1}, the following
distribution will play a role. Let $\Bsi\in\cD'$ be the distribution
which has the Fourier coefficients
\bqn
\Bsi_n &=& {\rm sign}(n).
\eqn
Moreover, let $\si\in\Cit{\{1\}}$ be the function
\bqn
\si(t) &=& \frac{1+t}{1-t}.
\eqn
Remark that both the distribution $\Bsi$ and the function $\si$ are odd.

\begin{proposition}\label{p1.2}
The distribution $\Bsi$ is in $\cD'(\{1\})$ and has the smooth part $\si$.
\end{proposition}
\begin{proof}
For $f\in\Cti{\{1\}}$ we can write $f(t)=(1-t)h(t)$ where $h\in\cD$.
Then
\bqn
[f\Bsi]_n &=&
\sum_{k>0}f_{n-k}-\sum_{k<0}f_{n-k}\nn\\
&=&
\sum_{k>0}(h_{n-k}-h_{n-k-1})-\sum_{k<0}(h_{n-k}-h_{n-k-1})\nn\\
&=&
h_{n-1}+h_{n}=
[h(t)(1+t)]_n=
[f\si]_n.\nn
\eqn
Note that $h_{k}$ converges to zero sufficiently fast.
This completes the proof.
\end{proof}

Given a distribution $\Ba\in\cD'$ with Fourier coefficients
$\{\Ba_n\}_{n=-\iy}^\iy$, define the $N\times N$ Toeplitz and Hankel
matrices by
\begin{equation}\label{f2.10}
T_N(\Ba) = \big(\Ba_{j-k}\big)_{j,k=0}^{N-1},\qquad\quad
H_N(\Ba) = \big(\Ba_{j+k+1}\big)_{j,k=0}^{N-1}.
\end{equation}
This definition is in accordance with (\ref{f1.2}).
Moreover define
\bqn\label{f2.11}
M_N(\Ba) = T_N(\Ba) + H_N(\Ba).
\eqn
Now we are ready to give the desired reformulation of Theorem \ref{t2.1}.

\begin{theorem}\label{t2.3}
Let $K$ be a compact and symmetric subset of $\T\setminus\{1\}$, and assume
that $\Ba\in\cD'(K)$ is an even distribution. Then
\bqn
\det T_{2N}(\Ba\Bsi) &=& (\det M_N(\Ba))^2.
\eqn
In the above $\Bc=\Ba\Bsi\in\cD'(K\cup\{1\})$ is an odd distribution.
\end{theorem}
\begin{proof}
Since $K$ and $\{1\}$ are disjoint compact sets, the distribution
$\Bc=\Ba\Bsi$ is well
defined. In the definition $\Bc=(\si f_{\si})\Ba+(af_a)\Bsi$ we may assume
without loss of generality that $f_\si$ and $f_a$ are even functions.
{}From this it follows easily that $\si f_{\si}$ is odd and $af_a$ is even;
thus both $(\si f_{\si})\Ba$ and $(af_a)\Bsi$ are odd. Hence
$\Bc$ is odd.

Next write $f_\si(t)=g(t)(1-t)(1-t\iv)$. Then
$(\si f_{\si})(t)=(1+t)(1-t\iv)g(t)=(t-t\iv)g(t)$. We obtain that
\bqn
[(\si f_{\si})\Ba]_n &=&
[((t-t\iv)g)\Ba]_n =
[g\Ba]_{n-1}-[g\Ba]_{n+1}.\nn
\eqn
Moreover, using the fact that $f_\si(t)=g(t)(1-t)(1-t\iv)$ and
keeping track of the cancellation, it follows
\bqn
\sum_{k=-n+1}^n [f_\si\Ba]_k
&=&
\sum_{k=-n+1}^n \Big(-[g\Ba]_{k-1}+2[g\Ba]_{k}-[g\Ba]_{k+1}\Big)
\nn\\
&=&
-[g\Ba]_{n}+[g\Ba]_{n-1}+[g\Ba]_{-n}-[g\Ba]_{-n-1}
\nn\\
&=&
[g\Ba]_{n-1}-[g\Ba]_{n+1}.\nn
\eqn
Here we have also used that $g\Ba$ is even.
{}From these two identities we obtain
\bqn\label{f1.14}
[(\si f_{\si})\Ba]_n &=&
\sum_{k=-n+1}^n [f_\si\Ba]_k.
\eqn
On the other hand, since $af_a=\Ba f_a$ is even,
\bqn\label{f1.15}
[(af_a)\Bsi)]_n &=&
\sum_{k>0}[af_a]_{n-k}-\sum_{k<0}[af_a]_{n-k}
=
\sum_{k=-n+1}^n[af_a]_{k}.
\eqn
Combining (\ref{f1.14}) and (\ref{f1.15}) yields
\bqn
\Bc_n &=& \sum_{k=-n+1}^n \Ba_k.\nn
\eqn
Together with Theorem \ref{t2.1} this completes the proof.
\end{proof}

Finally we will need the following result.

\begin{proposition}\label{p1.4}
For each $N\ge1$ we have $\det T_{2N}(\Bsi)=1$ and
\bqn\label{f1.16}
T_{2N}\iv(\Bsi) &=&
T_{2N}(\Bnu)
\eqn
where $\Bnu$ is the distribution with Fourier coefficients
$\Bnu_n={\rm sign}(n)(-1)^n$.
\end{proposition}
\begin{proof}
By Theorem \ref{t2.3} we have $\det T_{2N}(\Bsi)=(\det M_N(1))^2$,
where obviously $M_N(1)=I_N$. The formula for the inverse of
$T_{2N}(\Bsi)$ can be easily checked.
\end{proof}


\section{A limit theorem for determinants of skewsymmetric Toeplitz
matrices}
\label{sec:4}

Theorem \ref{t2.3} reduces the computation of the asymptotics
of $\det M_N(\Ba)$, which is what we are interested in, for certain
distributions $\Ba$, to the computation of the
asymptotics of $\det T_{2N}(\Bsi\Ba)$.

At first glance one might think that the asymptotics of $\det
T_{2N}(\Bsi\Ba)$
could be obtained from the predictions of the Fisher-Hartwig conjecture,
which was proved in \cite{E1}. Unfortunately, since $T_{2N}(\Bsi\Ba)$
is a skewsymmetric Toeplitz matrix, the Toeplitz determinant belongs to
those classes of functions where the Fisher-Hartwig conjecture breaks down.
It might be
the case that the asymptotic behavior fits with the
still unproved generalized conjecture
\cite{BT,E2}. However,
distributions of the kind $\Bsi\Ba$ appear here probably for the first time
in connection with Toeplitz determinants and since previous results do
not include this setting, previous techniques must be modified.

To that end, the purpose of this section is to prove that under certain
assumptions on the distribution $\Ba$ the expression
\bqn\label{f2.1}
\frac{\det T_{2N}(\Ba\Bsi)}{\det T_{2N}(\Ba)}
\eqn
converges to a certain (explicitly given) nonzero constant as $N\to\iy$.
Although we cannot rely on the main results of \cite{E1} (see also
\cite{E2}),
i.e., the Fisher-Hartwig conjecture, we will very heavily rely on the
machinery and several auxiliary results established in \cite{E1}.

Let us proceed with recalling the necessary definitions.
For $\mu\in\R$ let $\el_\mu$ stand for the Hilbert space of all
sequences $\{x_n\}_{n=0}^\iy$ of complex numbers for which
\bqn
\left\|\{x_n\}_{n=0}^\iy\right\|_{\mu} &:=&
\left(\sum_{n=0}^\iy(1+n)^{2\mu}|x_n|^{2}\right)^{1/2}
<\iy.
\eqn
For $\mu_1>\mu_2$ the space $\el_{\mu_1}$ is continuously and
densely embedded in $\el_{\mu_2}$.

The Toeplitz and the Hankel operator generated by $\Ba\in\cD'$
are the one-sided infinite matrices
\bqn
T(\Ba) = \left(\Ba_{j-k}\right)_{j,k=0}^\iy,\qquad
H(\Ba) = \left(\Ba_{j+k+1}\right)_{j,k=0}^\iy,
\eqn
where $\Ba_{n}$ are the Fourier coefficients of the distribution $\Ba$.
For each $\Ba\in\cD'$ there exist a (sufficiently large) $\mu_1$ and
a (sufficiently small) $\mu_2$ such that the matrices
$T(\Ba)$ and $H(\Ba)$ represent linear bounded operators acting from
$\el_{\mu_1}$ into $\el_{\mu_2}$.

The situation of the boundedness of Toeplitz and Hankel operators
generated by functions in $\cD$ was established in the following lemma
taken from \cite[Sect.~6.2]{E1}.

\begin{lemma}\label{l2.1}
For each $\mu,\mu_1,\mu_2\in\R$ and $a\in\cD$, the operator $T(a)$ is
bounded on $\cL(\el_\mu,\el_\mu)$ and the operator $H(a)$ is
bounded on $\cL(\el_{\mu_1},\el_{\mu_2})$.
\end{lemma}

We define the following finite rank operators acting on $\el_\mu$:
\bqn
P_N&:&(x_0,x_1,x_2,\dots)\mapsto
(x_0,x_1,\dots,x_{N-2},x_{N-1},0,0,\dots),\\
W_N&:&(x_0,x_1,x_2,\dots)\mapsto
(x_{N-1},x_{N-2},\dots,x_{1},x_{0},0,0,\dots).
\eqn
Obviously, $P_N^2=W_N^2=P_N$ and $W_NP_N=P_NW_N=W_N$.
If we consider the matrix $T_N(\Ba)$ as acting on the image
of the projection $P_N$ in the space $\el_{\mu}$, then
$T_N(\Ba)=P_NT(\Ba)P_N$.
Moreover,
\bqn\label{f3.6n}
W_NT_N(\Ba)W_N &=& T_N(\tilde{\Ba}).
\eqn
Recall that $\tilde{\Ba}$ is the distribution with the Fourier coefficients
$\tilde{\Ba}_{n}=\Ba_{-n}$.

For our purposes we need to single out two additional classes of
distributions.
Let $\Dp$ ($\Dm$, resp.) stand for the set of all distributions
$\Ba\in\cD'$ for which $\Ba_n=0$ for all $n<0$ ($n>0$, resp.). These two
sets form
commutative algebras with a unit element $e(t)\equiv1$. For $\Ba,\Bb\in\Dp$,
the product $\Bc=\Ba\Bb$ is defined by stipulating $\Bc_n=0$ for $n<0$ and
\bqn
\Bc_{n} &=& \sum_{k=0}^n \Ba_{n-k}\Bb_k\qquad
\mbox{ for } n\ge0.
\eqn
For $\Ba,\Bb\in\Dm$, the product $\Bc=\Ba\Bb$ is defined
by stipulating $\Bc_n=0$ for $n>0$ and
\bqn
\Bc_{n} &=& \sum_{k=n}^0 \Ba_{n-k}\Bb_k\qquad
\mbox{ for } n\le0.
\eqn
This definition of a multiplication is compatible with that of (\ref{f1.10})
whenever both are defined. Let $\cG\Dpm$ stand for the group of all
invertible distributions in $\Dpm$. Moreover, we put
$\Dpm(K)=\Dpm\cap\cD'(K)$ and let $\cG\Dpm(K)$ stand for the group of all
invertible elements in $\Dpm(K)$. One can show that if
$\Ba\in\cG\Dpm(K)$ has the smooth part $a$, then $a$ is an invertible
element of $\Cit{K}$.

There are some obvious relations between the distributions
$\Ba$ and $\tilde{\Ba}$. For instance, if $\Ba\in\cG\Dp(K)$, then
$\tilde{\Ba}\in\cG\Dm(\wt{K})$. Moreover, if $\Ba$ has the smooth part
$a$, then $\tilde{\Ba}$ has the smooth part $\tilde{a}$.

Let $H_1$ and $H_2$ be Hilbert spaces. We consider sequences
$\{C_N\}_{N=1}^\iy$ the elements of
which are well defined linear bounded operators (or matrices)
$C_N\,:\,H_1\to H_2$ for all sufficiently large $N$. Let
$\cO(\ro)$ with $\ro\in\R$ stand for the set
of all such sequences for which
\bqn\label{f3.10}
\|C_N\|_{\cL(H_1,H_2)} &=& O(N^\ro)\qquad
\mbox{ as }N\to\iy.
\eqn
The dependence of $\cO(\ro)$ on $H_1$ and $H_2$ will not
be displayed in the notation. We also use the notation
$\cO(\ro)$ in order to denote any sequence
of this type. In this sense,
$C_N=C+\cO(\ro)$ means that $\{C_N-C\}_{N=1}^\iy\in\cO(\ro)$.

Now let $H_1$, $H_2$, $\wt{H}_1$ and $\wt{H}_2$ be Hilbert spaces
and $\ros\in\R$. We denote by $\cO(\ros)$
the set of all sequences $\{C_N\}_{N=1}^\iy$ of $2\times 2$
block operators for which
\bqn
C_N&=&\left(\ba{cc}
\cO(\ro_0) & \cO(\ro_1) \\ \cO(\ro_2) & \cO(\ro_0)
\ea\right)\;:\;H_1\oplus\wt{H}_1\to H_2\oplus\wt{H}_2.
\eqn

We will also use the notations $\cO_2(\ro)$ and $\cO_1(\ro)$.
The only difference in comparison with $\cO(\ro)$ is that we consider the
convergence in (\ref{f3.10}) in the Hilbert-Schmidt and in the trace
class norm, respectively. Likewise, we will use the notations
$\cO_2(\ros)$ and $\cO_1(\ros)$.

Given $\Ba\in\cD'$, assume  that $T_N(\Ba)$ is invertible
for all sufficiently large $N$, and introduce the following
sequences of operators of $2\times 2$ block form:
\bqn\label{RR.1}
R_N(\Ba) &=& \left(\ba{cc} T_N\iv(\Ba) & T_N\iv(\Ba)W_N \\
T_N\iv(\tilde{\Ba})W_N & T_N\iv(\tilde{\Ba}) \ea \right),
\\[2ex]
RH_N(\Ba) &=& \left(\ba{c@{\quad}c}
T_N\iv(\Ba)P_NH(\Ba) & T_N\iv(\Ba)W_NH(\tilde{\Ba})\\
T_N\iv(\tilde{\Ba})W_NH(\Ba) & T_N\iv(\tilde{\Ba})P_NH(\tilde{\Ba})\ea
\right),\\[1ex]
HR_N(\Ba) &=& \left(\ba{c@{\quad}c}
H(\tilde{\Ba})P_NT_N\iv(\Ba) & H(\tilde{\Ba})W_NT_N\iv(\tilde{\Ba})\\
H(\Ba)W_NT_N\iv(\Ba) & H(\Ba)P_NT_N\iv(\tilde{\Ba})\ea
\right),\\[1ex]\qquad
HRH_N(\Ba)&=&
\left(\ba{c@{\quad}c}
H(\tilde{\Ba})P_NT_N\iv(\Ba)P_NH(\Ba) &
H(\tilde{\Ba})P_NT_N\iv(\Ba)W_NH(\tilde{\Ba})\\
H(\Ba)P_NT_N\iv(\tilde{\Ba})W_NH(\Ba) &
H(\Ba)P_NT_N\iv(\tilde{\Ba})P_NH(\tilde{\Ba})\ea
\right)\nn\\[1ex]
&&\;-\;\left(\ba{cc}T(\tilde{\Ba}) & H(\tilde{\Ba}\chi_{-N}) \\
H(\Ba\chi_{-N}) & T(\Ba) \ea\right).
\eqn
Here, as before,  $\chi_{-N}(t)=t^{-N}$, $t\in\T$.
These sequences of operators are considered
from $\el_{{\mu_1}}\oplus\el_{{\mu_1}}$
into $\el_{{\mu_2}}\oplus\el_{{\mu_2}}$
with ${\mu_1}$ sufficiently large and ${\mu_2}$ sufficiently small,
which ensures the boundedness of the operators.

Now we are prepared to define the notion of $\cR$-convergence,
which has been introduced in \cite{E1,E2}.
Let $\Ba\in\cD'$, $\Ba_+\in\cG\Dp$ and $\Ba_-\in\cG\Dm$ be distributions,
and let
$\ros\in\R$. We say that the distribution $\Ba$
{\em effects $\cR$-convergence with respect to $[\Ba_+,\Ba_-]$ and $(\ros)$}
if there exist ${\mu_1}\geq0$ and ${\mu_2}\leq0$ such that
\bqn\label{f112.1}
R_N(\Ba)&=&
\diag\Big(T(\Ba_+\iv)T(\Ba_-\iv),\,T(\tilde{\Ba}_-\iv)T(\tilde{\Ba}_+\iv)\Big)
+\cO(\ros),\qquad\\[.5ex]\label{f112.2}
RH_N(\Ba)&=&
\diag\Big(T(\Ba_+\iv)H(\Ba_+),\,T(\tilde{\Ba}_-\iv)H(\tilde{\Ba}_-)\Big)
+\cO(\ros),\\[.5ex]\label{f112.3}
HR_N(\Ba)&=&
\diag\Big(H(\tilde{\Ba}_-)T(\Ba_-\iv),\,H(\Ba_+)T(\tilde{\Ba}_+\iv)\Big)
+\cO(\ros),\\[.5ex]\label{f112.4}
HRH_N(\Ba)&=&\mbox{}-
\diag\Big(T(\tilde{\Ba}_-)T(\tilde{\Ba}_+),\,T(\Ba_+)T(\Ba_-)\Big)
+\cO(\ros),
\eqn
where these sequences are considered from
$\el_{{\mu_1}}\oplus\el_{{\mu_1}}$ into
$\el_{{\mu_2}}\oplus\el_{{\mu_2}}$.

We will need the concept of $\cR$-convergence for a distribution $\Ba$ which
is
even. This particular case gives rise to some simplifications.
For any distribution $\Ba\in\cD'$, the following statement are equivalent:
\begin{itemize}
\item[(1)]
$\Ba$ effects $\cR$-convergence with respect to
$[\Ba_+,\Ba_-]$ and $(\ros)$;
\item[(2)]
$\tilde{\Ba}$ effects $\cR$-convergence with respect to
$[\tilde{\Ba}_-,\tilde{\Ba}_+]$ and $(\ro_0,\ro_2,\ro_1)$.
\end{itemize}
In fact, in order to prove this equivalence, one has only to pass
to the transpose in equations (\ref{f112.1})--(\ref{f112.4}).
Hence for an even distribution $\Ba$ we can replace $\ro_1$ and $\ro_2$
by $\min\{\ro_1,\ro_2\}$, i.e., we may assume that $\ro_1=\ro_2=:\ro$.
Moreover, since the distributions $\Ba_+$ and $\Ba_-$ are uniquely
determined
up to a nonzero multiplicative constant, we may assume without loss of
generality that $\Ba_-=\tilde{\Ba}_+$.

This last remark gives some motivation for the assumptions in
the following theorem. In this theorem we establish the
asymptotic formula for (\ref{f2.1}).

\begin{theorem}[Limit Theorem]\label{t2.2}
Let $K$ be a symmetric and compact subset of $\T\setminus\{1,-1\}$.
Moreover, assume that
\begin{itemize}
\item[(i)]
$\Ba\in\cD'(K)$ is an {\em even} distribution with the smooth part
$a\in\Cit{K}$;
\item[(ii)]
$\Ba_+\in\cG\Dp(K)$ is a distribution with the smooth part $a_+\in\Cit{K}$;
\item[(iii)]
$a(t)=a_+(t)\tilde{a}_+(t)$ for all $t\in\T\setminus K$;
\item[(iv)]
$\ro_0<0$ and $\ro<0$;
\item[(v)]
$\Ba$ effects $\cR$-convergence with respect to
$[\Ba_+,\tilde{\Ba}_+]$ and $(\ro_0,\ro,\ro)$.
\end{itemize}
Then
\bqn
\frac{\det T_{2N}(\Ba\Bsi)}{\det T_{2N}(\Ba)}
&=& \frac{a_+(1)}{a_+(-1)} +O(N^{\max\{\ro_0,\ro\}})\qquad
\mbox{ as }N\to\iy.
\eqn
\end{theorem}

The rest of this section is devoted to the proof of this theorem.
Once again we need to quote some auxiliary results.

Let us introduce the notation
\bqn
X^\mi_\mii=\el_\mi\oplus\el_\mii,\qquad
X^\mi=\el_\mi,\qquad
X_\mii=\el_\mii.
\eqn
This notation is convenient in the sense that it reflects the condition
that $\mi$ is a sufficiently large and $\mii$ a sufficiently small real
number, which we will encounter in what follows.

The following proposition is taken from \cite[Proposition~8.1]{E1}
with a slight change of notation. As before, we will denote
a distribution by a bold letter and its smooth part by the same non-bold
letter without mentioning it explicitly.

\begin{proposition}\label{p2.3}
Let $K$ be a compact subset of $\T$, and let $\Ba\in\cD'(K)$ and
$f\in\Cti{K}$.
Then for all sufficiently large $\mi\ge0$ and all sufficiently small
$\mii\le0$ the linear operators
\bqn
H_1(\Ba,f) &=&
\Big(H(\Ba),\; H(af)-T(\Ba)H(f)\;\Big),\\
H_2(f,\Ba) &=&
\left(\ba{c} H(\tilde{f}\tilde{a})-H(\tilde{f})T(\Ba)\\
H(\tilde{\Ba})\ea\right)
\eqn
are bounded on the spaces
\bqn
H_1(\Ba,f):X^\mi_\mii\to X_\mii,\qquad
H_2(f,\Ba):X^\mi\to X^\mi_\mii.
\eqn
\end{proposition}
In order to provide some meaning to these operators we remark that
with appropriately chosen $\mi$ and $\mii$, these operators can be
embedded into the spaces
\bqn
H_1(\Ba,f):X^\mi_\mi\to X_\mii,\qquad
H_2(f,\Ba):X^\mi\to X^\mii_\mii.
\eqn
In this case, these operators can be written as a product:
\bqn
H_1(\Ba,f) &=& H(\Ba)\Big(I,\; T(\tilde{f})\;\Big),\\
H_2(f,\Ba) &=&
\left(\ba{c} T(\tilde{f})\\ I\ea\right)H(\tilde{\Ba}).
\label{f2.44}
\eqn

The next result is taken from \cite[Corollary~8.8]{E1}.

\begin{proposition}\label{p2.4}
Let  $K_1$ and $K_2$ be disjoint and compact subsets of $\T$, and let
$\Ba_i\in\cD'(K_i)$ and $f_i\in\Cti{K_i}$ $(i=1,2)$ be such that
$f_1+f_2=1$.
Then
\bqn
T_N(\Ba_1\Ba_2) &=&
T_N(\Ba_1)T_N(\Ba_2) +
P_NH_1(\Ba_1,f_1)H_2(f_2,\Ba_2)P_N\nn\\[1ex]
&&\mbox{}+
W_NH_1(\tilde{\Ba}_1,\tilde{f}_1)H_2(\tilde{f}_2,\tilde{\Ba}_2)W_N.
\label{f2.45}
\eqn
\end{proposition}
We remark that the linear operators occurring in (\ref{f2.45})
are bounded on appropriately chosen spaces. Moreover, (\ref{f2.45})
represents
a generalization of the well-known identity
\bqn
T_{N}(a_1a_2) &=& T_N(a_1)T_N(a_2) +
P_{N}H(a_1)H(\tilde{a}_2)P_N+
W_{N}H(\tilde{a}_1)H(a_2)W_N\nn
\eqn
due to Widom \cite{W1}, which holds for functions $a_1,a_2\in L^\iy(\T)$.

Next, we consider the functions $\xi_1(t)=1-t\iv$ and $\xi_{-1}(t)=1+t\iv$.
These functions can be identified with distributions $\Bxi_1\in\Dm$
and $\Bxi_{-1}\in\Dm$, respectively, in the sense of (\ref{f2.4}).
Obviously, both $\Bxi_1$ and $\Bxi_{-1}$ belong to $\cG\Dm$.
In fact, the inverse distributions $\Bxi\iv_{1}$ and $\Bxi\iv_{-1}$
are given as follows by their Fourier coefficients:
\bqn
\left[\Bxi\iv_{1}\right]_n &=& \left\{
\ba{ll} 0&\mbox{ if }n>0\\ 1&\mbox{ if } n\le0,\ea\right.\\
\left[\Bxi\iv_{-1}\right]_n &=& \left\{
\ba{ll} 0&\mbox{ if }n>0\\ (-1)^n&\mbox{ if } n\le0.\ea\right.
\eqn
\begin{proposition}\label{p2.5}
The following statements are true:
\begin{itemize}
\item[(a)]
$\Bxi_1\in\cG\Dm(\{1\})$, and $\Bxi_1^{\pm1}$ has the smooth part
$\xi_1^{\pm1}$;
\item[(b)]
$\Bxi_{-1}\in\cG\Dm(\{-1\})$, and $\Bxi_{-1}^{\pm1}$ has the smooth part
$\xi_{-1}^{\pm1}$.
\end{itemize}
\end{proposition}
\begin{proof}
Since $\Bxi_1=\xi_1$ in the sense of (\ref{f2.4}), we even have
$\Bxi_1\in\cD'(\emptyset)$ and $\Bxi_1$ has the smooth part $\xi_1$.
It remains to show that $\Bxi_1\iv$ is contained in $\cD'(\{1\})$
and has the smooth part $\xi_1\iv$. Indeed, let $f\in\Cti{\{1\}}$.
Then we can write $f=\xi_1g$ where $g\in C^\iy(\T)$.
It can be checked easily that $\Bxi_1\iv\xi_1=1$.
Hence $\Bxi_1\iv f=g=\xi_1\iv f$. This completes the proof of part (a).
Part (b) can be proved analogously.
\end{proof}

The following corollary, in which we define another distribution $\Bh$,
is a simple consequence of the previous proposition.
\begin{corollary}\label{c2.6}
The distribution $\Bh=\Bxi_1\iv\Bxi_{-1}$ is contained in $\cG\Dm(\{-1,1\})$
and has the smooth part $h=\xi_1\iv\xi_{-1}$. The inverse distribution
$\Bh\iv$ equals $\Bxi_1\Bxi_{-1}\iv$ and has the smooth part
$h\iv=\xi_1\xi_{-1}\iv$.
\end{corollary}

It is easy to compute the Fourier coefficients of $\Bh$ and $\Bh\iv$:
\bqn
[\Bh]_n &=&
\left\{\ba{ll} 0&\mbox{ if }n>0\\
1&\mbox{ if } n=0\\ 2&\mbox{ if }n<0,\ea\right.
\\ \label{f2.50x}
[\Bh\iv]_n &=&
\left\{\ba{ll} 0& \mbox{ if }n>0\\
1&\mbox{ if } n=0\\ 2(-1)^n&\mbox{ if }n<0.\ea\right.
\eqn

There is (for our purposes) an important relation between operators
containing the distributions $\Bsi$ and $\Bh$,
which is given in the following proposition.

\begin{proposition}\label{p2.7}
Let $f\in\Cti{\{-1,1\}}$. Then
\bqn
&&
H_2(f,\Bsi)P_{2N}T_{2N}\iv(\Bsi) =
H_2(f,\Bsi)W_{2N}T_{2N}\iv(\Bsi) =
\nn\\
&& =
\frac{1}{2}H_2(f,\Bh)P_{2N}T_{2N}\iv(\Bh) =
-\frac{1}{2}H_2(f,\Bh)P_{2N}T_{2N}\iv(\Bh)W_{2N}.
\label{f3.32}
\eqn
\end{proposition}
\begin{proof}
We first prove that
\bqn
&&
H(\tilde{\Bsi})P_{2N}T_{2N}\iv(\Bsi) =
H(\tilde{\Bsi})W_{2N}T_{2N}\iv(\Bsi) =
\nn\\
&& =
\frac{1}{2} H(\tilde{\Bh})P_{2N}T_{2N}\iv(\Bh) =
-\frac{1}{2} H(\tilde{\Bh})P_{2N}T_{2N}\iv(\Bh)W_{2N}.
\label{f2.49}
\eqn
In this identity we do not need to worry about the boundedness on certain
spaces since both the left and right hand side is
an infinite Hankel matrix times a finite rank matrix.
Let $x=(1,1,\dots)^T$ denote an infinite column vector and
$x_{2N}=(1,\dots,1)^T$ a finite column vector of size $2N$.
Then
\begin{equation}
H(\tilde{\Bsi})P_{2N} = H(\tilde{\Bsi})W_{2N} = -xx_{2N}^T,\qquad
H(\tilde{\Bh})P_{2N} = 2 xx_{2N}^T.
\end{equation}
A moments thought shows that (\ref{f2.49}) is proved
as soon as
\begin{equation}
-x_{2N}^TT_{2N}\iv(\Bsi) = x_{2N}^TT_{2N}\iv(\Bh)
= -x_{2N}^TT_{2N}\iv(\Bh)W_{2N}
\end{equation}
is established. However, this is just a straightforward calculation.
We have to observe that $T_{2N}\iv(\Bh)=T_{2N}(\Bh\iv)$ with the
Fourier coefficients given by (\ref{f2.50x}) and moreover that
$T_{2N}\iv(\Bsi)=T_{2N}(\Bnu)$ by Proposition \ref{p1.4}.

Having proved (\ref{f2.49}), we take into account
the identity (\ref{f2.44}) and formula (\ref{f3.32}) follows by a
density argument of the Hilbert spaces under consideration.
\end{proof}

In \cite[Formula (8.64)]{E1}, the following $2\times2$ block operators
acting on $X^\mi_\mii\oplus X^\mi_\mii$ with sufficiently large $\mi$
and small $\mii$ were defined:
\bqn
\lefteqn{HSH_N(f_2,\Ba_2,\Ba_1,f_1)}\hspace{0ex}
\label{f3.36}\\
&=&
\left(\ba{c} H_2(f_2,\Ba_2)P_{N} \\ H_2(\tilde{f}_2,\tilde{\Ba}_2)W_{N}
\ea\right)T_N\iv(\Ba_2)T_N\iv(\Ba_1)
\Big(P_NH_1(\Ba_1,f_1),\;\;W_NH_1(\tilde{\Ba}_1,\tilde{f}_1)\;
\Big).\nn
\eqn
Moreover, it has been shown \cite[Formula (8.139)]{E1} that
\bqn\label{f2.HSH}
HSH_N(f_2,\Ba_2,\Ba_1,f_1) &=&
HZ_N(f_2,\Ba_2,f_2)YH_N(f_1,\Ba_1,f_1),
\eqn
where $HZ_N(\dots)$ and $YH_N(\dots)$ are linear bounded
$2\times2$ block operators acting on $X^\mi_\mii\oplus X^\mi_\mii$,
which we are not going to define here.
The important result concerning these operators is the following asymptotic
formula, which is taken from \cite[Proposition~10.4 and
Proposition~10.5]{E1}.

\begin{proposition}\label{p2.8}
Let $K$ be a compact subset of $\T$, $f\in\Cti{K}$, and assume that
$\Ba\in\cD'(K)$, $\Ba_\pm\in\cG\Dpm(K)$ such that $a=a_+a_-$ holds
for their smooth parts.
If $\Ba$ effects $\cR$-convergence with respect to
$[\Ba_+,\Ba_-]$ and $(\ros)\in\R^3$, then
\bqn
YH_N(f,\Ba,f) &=& \diag\left(TH(f,\Ba_+,f),\;TH(\tilde{f},\tilde{\Ba}_-,
\tilde{f})\right)+\cO_2(\ros),\nn\\
HZ_N(f,\Ba,f) &=& \diag\left(HT(f,\Ba_-,f),\;HT(\tilde{f},\tilde{\Ba}_+,
\tilde{f})\right)+\cO_2(\ros).\nn
\eqn
\end{proposition}

In the previous proposition linear operators $TH(\dots)$ and $HT(\dots)$
appear, which were defined in \cite[Formulas (9.11) and (9.12)]{E1}.
These operators are Hilbert-Schmidt operators
on the space $X^\mi_\mii$ for all
sufficiently large $\mi\ge0$ and all sufficiently small $\mii\le0$
(see \cite[Proposition~9.2]{E1}).

Moreover, the asymptotic operator relation stated in the previous
proposition has to be understood in the way that the operators
act on the space $X^\mi_\mii\oplus X^\mi_\mii$,
where $\mi\ge0$ is fixed and sufficiently large and $\mii\le0$ is fixed and
sufficiently small.

Now we are prepared to give the proof of Theorem \ref{t2.2}.

\begin{proofof}{Theorem \ref{t2.2}}
We start from Proposition \ref{p2.4} with $\Ba_1=\Ba$, $K_1=K$,
$\Ba_2=\Bsi$, $K_2=\{-1,1\}$. Since $K$ and $\{1,-1\}$ are
symmetric subsets of $K$, we can assume without loss of generality that
$f_1$ and $f_2$ are even functions. Then
\bqn
T_{2N}(\Ba\Bsi) &=&  T_{2N}(\Ba)T_{2N}(\Bsi)+
P_{2N}H_1(\Ba,f_1)H_2(f_2,\Bsi)P_{2N}\nn\\
&&\mbox{}-
W_{2N}H_1(\Ba,f_1)H_2(f_2,\Bsi)W_{2N},
\eqn
where we have also used that $\Ba$ is even and $\Bsi$ is odd.
By Proposition \ref{p1.4} the inverses of $T_{2N}(\Bsi)$ exist for all
$N$. Since the distribution $\Ba$ effects $\cR$-convergence, the inverses of
$T_{2N}(\Ba)$ exist for all sufficiently large $N$. Hence
\bqn
T_{2N}\iv(\Ba)T_{2N}(\Ba\Bsi)T_{2N}\iv(\Bsi)
&=& P_{2N}+
T_{2N}\iv(\Ba)P_{2N}H_1(\Ba,f_1)H_2(f_2,\Bsi)P_{2N}T_{2N}\iv(\Bsi)
\nn\\&&\mbox{}-
T_{2N}\iv(\Ba)W_{2N}H_1(\Ba,f_1)H_2(f_2,\Bsi)W_{2N}T_{2N}\iv(\Bsi).\nn
\eqn
{}From Proposition \ref{p2.7} it now follows that
\bqn
T_{2N}\iv(\Ba)T_{2N}(\Ba\Bsi)T_{2N}\iv(\Bsi)
&=&
P_{2N}+\frac{1}{2}
T_{2N}\iv(\Ba)P_{2N}H_1(\Ba,f_1)H_2(f_2,\Bh)P_{2N}T_{2N}\iv(\Bh)
\nn\\&&\mbox{}-\frac{1}{2}
T_{2N}\iv(\Ba)W_{2N}H_1(\Ba,f_1)H_2(f_2,\Bh)P_{2N}T_{2N}\iv(\Bh).\nn
\eqn
Taking determinants, observing that $\det T_{2N}(\Bsi)=1$ by Proposition
\ref{p1.4} and using the formula $\det(I+AB)=\det(I+BA)$ for determinants,
we obtain
\bqn
\lefteqn{\frac{\det T_{2N}(\Ba\Bsi)}{\det T_{2N}(\Ba)}}
\nn\\
&=&
\det\left(
P_{2N}+\frac{1}{2}
T_{2N}\iv(\Ba)(P_{2N}-W_{2N})H_1(\Ba,f_1)
H_2(f_2,\Bh)P_{2N}T_{2N}\iv(\Bh)
\right)\nn\\
&=&
\det\left(
I+\frac{1}{2}H_2(f_2,\Bh)P_{2N}T_{2N}\iv(\Bh)
T_{2N}\iv(\Ba)(P_{2N}-W_{2N})H_1(\Ba,f_1)
\right)\nn\\
&=&
\det\left(
I+\frac{1}{2}H_2(f_2,\Bh)P_{2N}T_{2N}\iv(\Bh)(P_{2N}-W_{2N})
T_{2N}\iv(\Ba)P_{2N}H_1(\Ba,f_1)
\right).\nn
\eqn
Here we have also used formula (\ref{f3.6n}) and the fact that $\Ba$
is even. Again from Proposition \ref{p2.7} it follows that
\bqn
\lefteqn{\frac{\det T_{2N}(\Ba\Bsi)}{\det T_{2N}(\Ba)}}
\nn\\
&=&
\det\left(
I+H_2(f_2,\Bh)P_{2N}T_{2N}\iv(\Bh)
T_{2N}\iv(\Ba)P_{2N}H_1(\Ba,f_1)\right)
\nn\\
&=&
\det\left(
I+\Big(I,\;0\Big)HSH_{2N}(f_2,\Bh,\Ba,f_1)\left(\ba{c}I\\0\ea\right)\right).
\label{f3.39}
\eqn
Note that in the last formula the $(1,1)$-block entry of
$HSH_{2N}(f_2,\Bh,\Ba,f_1)$ appears (see formula (\ref{f3.36})).

{}From the assumption on $\Ba$ and -- as concerns the distribution
$\Bh$ -- from Proposition \ref{p2.5}, Corollary \ref{c2.6} and
\cite[Theorem 13.1]{E1} (see also the remark made after
Theorem \ref{t4.3new} below) we know that
\begin{itemize}
\item[(1)]
$\Ba$ effects $\cR$-convergence with respect to $[\Ba_+,\tilde{\Ba}_+]$
and $(\ro_0,\ro,\ro)$;
\item[(2)]
$\Bh$ effects $\cR$-convergence with respect to $[1,\Bh]$
and $(\omega,0,\omega)$ for each $\omega\in\R$.
\end{itemize}
Hence, by Proposition \ref{p2.8} we can conclude that
\bqn
YH_{2N}(f_1,\Ba,f_1) &=&
\left(\ba{cc}TH(f_1,\Ba_+,f_1)+\cO_2(\ro_0) & \cO_2(\ro)\\
\cO_2(\ro)&TH(\tilde{f}_1,\tilde{\Ba}_+,\tilde{f}_1)+\cO_2(\ro_0)
\ea\right),\nn\\
HZ_{2N}(f_1,\Bh,f_2) &=&
\left(\ba{cc}HT(f_2,\Bh,f_2)+\cO_2(\omega) & \cO_2(0)\\
\cO_2(\omega)&HT(\tilde{f}_2,1,\tilde{f}_2)+\cO_2(\omega)
\ea\right).\nn
\eqn
This in connection with (\ref{f2.HSH}) yields
\bqn
\lefteqn{\Big(I,\;0\Big)HSH_{2N}(f_2,\Bh,\Ba,f_1)
\left(\ba{c}I\\0\ea\right)}\nn\\
&=& \Big(HT(f_2,\Bh,f_2)+\cO_2(\omega)\Big)\Big(
TH(f_1,\Ba_+,f_1)+\cO_2(\ro_0)\Big) + \cO_2(0)\cO_2(\ro).\nn
\eqn
Noting that the operators $HT(\dots)$ and $TH(\dots)$ are Hilbert-Schmidt
and choosing $\omega$ sufficiently small, this implies
\bqn
\lefteqn{\Big(I,\;0\Big)HSH_{2N}(f_2,\Bh,\Ba,f_1)
\left(\ba{c}I\\0\ea\right)}\nn\\
&=&
HT(f_2,\Bh,f_2)TH(f_1,\Ba_+,f_1)+\cO_1(\max\{\ro_0,\ro\}).
\label{f3.40}
\eqn
Formulas (\ref{f3.39}) and (\ref{f3.40}) give
\bqn
\frac{\det T_{2N}(\Ba\Bsi)}{\det T_{2N}(\Ba)} &=&
\det\left(I+HT(f_2,\Bh,f_2)TH(f_1,\Ba_+,f_1)\right)
+O(N^{\max\{\ro_0,\ro\}}).\nn
\eqn

It remains to show that the above operator determinant equals the constant
$a_+(1)/a_+(-1)$. From \cite[Proposition~9.10(b)]{E1} we obtain
\bqn
\det\left(I+HT(f_2,\Bh,f_2)TH(f_1,\Ba_+,f_1)\right)
&=& \lim_{r\to1-0} E(h_r\Ba_+,h_r\Bh)
\eqn
where $h_r\Ba_+$ and $h_r\Bh$ ($0\le r<1$) are the harmonic extensions of
the
distributions $\Ba_+$ and $\Bh$, i.e.,
\begin{equation}
(h_r\Ba_+)(t)=\sum_{n=0}^\iy r^nt^n[\Ba_+]_n,\qquad
(h_r\Bh)(t)=\sum_{n=0}^\iy r^nt^{-n}\Bh_{-n},
\end{equation}
and $E(\dots)$ is the constant defined by
\bqn
E(\phi_+,\phi_-) &=& \exp\left(\sum_{n=1}^\iy n[\log\phi_+]_n
[\log\phi_-]_{-n}\right).
\eqn
In this connection we remark that the harmonic extensions of distributions
contained in $\cG\Dp$ or in $\cG\Dm$ are always functions in $\cD=\Cti{}$,
which possess a continuous logarithm on $\T$ for each $0\le r<1$. Indeed,
if $\Bb\in\cG\Dpm$, then the harmonic extensions are multiplicative, and
consequently $(h_r\Bb)(h_r\Bb\iv)\equiv1$. The harmonic extensions depend
uniformly on $r$ and are constants for $r=0$. Hence the functions
$h_r\Bb$ are nonzero on all of $\T$ and have winding number zero.

In order to compute $E(h_r\Ba_+,h_r\Bh)$, observe first that
$$
h_r\Bh=(h_r\Bxi_1\iv)(h_r\Bxi_{-1})=
\left(1-\frac{r}{t}\right)\iv\left(1+\frac{r}{t}\right).
$$
Hence
\bqn
[\log(h_r\Bh)]_{-n}&=&\frac{r^n}{n}-\frac{(-r)^n}{n},\qquad n\ge1.
\eqn
We obtain that $E(h_r\Ba_+,h_r\Bh)$ is the exponential of
\bqn
\lefteqn{
\sum_{n=1}^\iy \left([\log h_r\Ba_+]_nr^n-[\log h_r\Ba_+]_n(-r)^n)\right)}
\hspace{8ex}\nn\\
&=&
(h_r(\log h_r\Ba_+))(1)-(h_r(\log h_ra_+))(-1).\nn
\eqn
Now notice that for $b\in\Cti{}\cap\Dp$ we have $h_re^b=\exp(h_rb)$.
With $b=\log h_r\Ba_+$ we obtain
$\exp(h_r(\log h_r\Ba_+))=h_r(h_r\Ba)=h_{r^2}\Ba_+$.
Hence
\bqn
E(h_r\Ba_+,h_r\Bh) &=&
\exp\left( (\log h_{r^2}\Ba_+)(1)-(\log h_{r^2}\Ba_+)(-1)\right).
\eqn
{}From \cite[Proposition~4.4]{E1} it follows that this converges to
$a_+(1)/a_+(-1)$. Thus the proof is complete.
\end{proofof}


\section{Asymptotics of the determinants of symmetric Toeplitz matrices
with Fisher-Hartwig distributions}
\label{sec:5}

In this section we recall the known results about the asymptotic behavior
of Toeplitz determinants and specialize them to the case of determinants
of symmetric Toeplitz matrices. Such an asymptotic formula is
provided by the Fisher-Hartwig conjecture, which -- in the so far
most general setting -- has been proved in \cite{E1} (see also \cite{E2}).
We also recall the some known results about the $\cR$-convergence
of certain classes of distributions, which are later on needed in order to
be
able to use Theorem \ref{t2.2}. The underlying classes of distributions
are defined and their properties are stated next.

For $\alpha,\beta,\gamma,\delta\in\C$ and $\theta_0\in(-\pi,\pi]$,
we introduce the functions
\bqn
\omega_{\alpha,\beta,\theta_0}(e^{i\theta}) &=&
(2-2\cos(\theta-\theta_{0}))^\alpha e^{i\beta(\theta-\theta_{0}-\pi)},
\qquad 0<\theta-\theta_{0}<2\pi.
\\
\eta_{\gamma,\theta_0}(e^{i\theta}) &=&
(1-e^{i(\theta-\theta_0)})^\gamma,
\\
\xi_{\delta,\theta_0}(e^{i\theta}) &=&
(1-e^{i(\theta_0-\theta)})^\delta.
\eqn
This implies
\bqn
\omega_{\alpha,\beta,\theta_0}(e^{i\theta}) &=&
\eta_{\alpha+\beta,\theta_0}(e^{i\theta})
\xi_{\alpha-\beta,\theta_0}(e^{i\theta}).
\eqn

For $\alpha,\beta,\gamma,\delta\in\C$ with
$2\alpha\notin\Zm:=\{-1,-2,-3,\dots\}$
and $\theta_0\in(-\pi,\pi]$, we introduce the distributions
$\Bom_{\alpha,\beta,\theta_0}$, $\Bet_{\gamma,\theta_0}$ and
$\Bxi_{\delta,\theta_0}$ in terms of their Fourier coefficients
\bqn
[\Bom_{\alpha,\beta,\theta_0}]_n
&=&
\frac{e^{in(\pi-\theta_0)}\Ga(1+2\alpha)}
{\Ga(1+\alpha+\beta-n)\Ga(1+\alpha-\beta+n)},
\qquad n\in\Z,
\\{}
[\Bet_{\gamma,\theta_0}]_n
&=&
\left\{\ba{ll}
e^{in(\pi-\theta_0)}{\gamma\choose n}
&\mbox{ if }n\ge0\\0&\mbox{ if }n<0,
\ea\right.
\\{}
[\Bxi_{\delta,\theta_0}]_n
&=&
\left\{\ba{ll}
0&\mbox{ if }n>0\\
e^{in(\pi-\theta_0)} {\delta\choose -n}
&\mbox{ if }n\le0.
\ea\right.
\eqn
It can be checked straightforwardly that if $2\alpha\not\in\Zm$, then
\bqn
\Bom_{\alpha,\alpha,\theta_0}=\Bet_{2\alpha,\theta_0},\qquad
\Bom_{\alpha,-\alpha,\theta_0}=\Bxi_{2\alpha,\theta_0}.
\eqn

In what follows let $\G\Cti{}$ stand for the set of all functions
in $\Cti{}$ which are nonzero on all of $\T$ and have winding number zero.
In other words, $\G\Cti{}$ is the set of all complex-valued functions
defined on $\T$ which possess a logarithm that belongs to $\Cti{}$.
Moreover, let $\Cti{\pm}$ stand for the set of all $f\in\Cti{}$ for which
$f_n=0$ for all $n<0$ ($n>0$, resp.). We denote by $G\Cti{\pm}$ the set
of all invertible functions in $\Cti{\pm}$.

A {\em distribution of Fisher-Hartwig type} is a distribution of the form
\bqn\label{f3.6}
\Bc &=& b\prod_{r\in M_0} \Bom_{\alpha_r,\beta_r,\theta_r}
\prod_{r\in M_+}\Bet_{\gamma_r,\theta_r}
\prod_{r\in M_-}\Bxi_{\delta_r,\theta_r},
\eqn
where
\begin{itemize}
\item[(i)]
$R\ge0$ and $\{1,\dots,R\}=M_0\cup M_+\cup M_-$
is a decomposition into disjoint subsets;
\item[(ii)]
$\theta_1,\dots,\theta_R\in(-\pi,\pi]$ are distinct numbers;
\item[(iii)]
$b\in \G\Cti{}$;
\item[(iv)]
$\alpha_r,\beta_r\in\C$ and $2\alpha_r\notin\Zm$ for all $r\in M_0$;
\item[(v)]
$\gamma_r\in\C$ for all $r\in M_+$ and $\delta_r\in\C$ for all $r\in M_-$.
\end{itemize}
The product (\ref{f3.6}) of the distributions has to be understood
in the sense of (\ref{f1.10}). To this distribution we associate the
function
\bqn\label{f3.7}
c(e^{i\theta}) &=& b(e^{i\theta})\prod_{r\in M_0}
\omega_{\alpha_r,\beta_r,\theta_r}(e^{i\theta})
\prod_{r\in M_+}\eta_{\gamma_r,\theta_r}(e^{i\theta})
\prod_{r\in M_-}\xi_{\delta_r,\theta_r}(e^{i\theta}).\qquad
\eqn
Such a function will be called a {\em function of Fisher-Hartwig type}.

The following result has been proved in \cite[Proposition 5.5]{E1}.
\begin{proposition}\label{p4.1new}
Let $\Bc$ be the distribution \refeq{f3.6} and $c$ be the
function \refeq{f3.7}. Put $K=\{e^{i\theta_r}:1\le r\le R\}$.
Then
\begin{itemize}
\item[(a)]
$\Bc\in\cD'(K)$ and $\Bc$ has the smooth part $c$;
\item[(b)]
if $M_0=M_-=\emptyset$ and $b\in G\Cti{+}$, then $\Bc\in\cG\Dp(K)$;
\item[(c)]
if $M_0=M_+=\emptyset$ and $b\in G\Cti{-}$, then $\Bc\in\cG\Dm(K)$.
\end{itemize}
Moreover, if $\Re\alpha>-1/2$, $\Re\gamma>-1$ and $\Re\delta>-1$, then
the distribution $\Bc$ can be identified with $c\in L^1(\T)$ in the sense
of \refeq{f2.4}.
\end{proposition}

In what follows, we agree on the following conventions. We say that
\begin{itemize}
\item[(i)]
$\Ba$ effects $\cR$-convergence w.r.t.\ $[\Ba_+,\Ba_-]$ and
$(-\iy,-\iy,-\iy)$ if and only if
for each $\ro\in\R$, $\Ba$ effects $\cR$-convergence w.r.t.\
$[\Ba_+,\Ba_-]$ and $(\ro,\ro,\ro)$;
\item[(ii)]
$\Ba$ effects $\cR$-convergence w.r.t.\ $[\Ba_+,a_-]$ and
$(-\iy,-\iy,\mu)$ if and only if
for each $\ro\in\R$, $\Ba$ effects $\cR$-convergence w.r.t.\
$[\Ba_+,\Ba_-]$ and $(\ro,\ro,\mu)$;
\item[(iii)]
$\Ba$ effects $\cR$-convergence w.r.t.\ $[\Ba_+,\Ba_-]$ and
$(-\iy,\mu,-\iy)$ if and only if
for each $\ro\in\R$, $\Ba$ effects $\cR$-convergence w.r.t.\
$[\Ba_+,\Ba_-]$ and $(\ro,\mu,\ro)$.
\end{itemize}
Finally, let $O(N^{-\iy})$ stand for a sequence of complex numbers which
is $O(N^\ro)$ for each $\ro\in\R$. A maximum taken oven
an empty set is considered to be $-\iy$.

Under certain conditions on the parameters, the asymptotic behavior
of the determinants $\det T_N(\Bc)$ with $\Bc$ given by (\ref{f3.6})
is described by the Fisher-Hartwig conjecture. The proof of this conjecture
together with the statement that such distributions effect
$\cR$-convergence was the main result of \cite{E1}
(see Section 13 therein or \cite[Section 6]{E2}).

Recall that $\Zm:=\{-1,-2,-3,\dots\}$. Moreover, given $b\in\G\Cti{}$,
let  $G[b]$ and $E[b]$ stand for the constants (\ref{f1.Gb}) and
(\ref{f1.Eb}),
and let $b_\pm\in G\Cti{\pm}$ stand for the functions (\ref{f1.bpm}).

\begin{theorem}\label{t4.3new}
Let $\theta_1,\dots,\theta_R\in(-\pi,\pi]$ be distinct numbers, $R\ge0$, put
$t_r=e^{i\theta_r}$, let $\{1,\dots,R\}=M_0\cup M_+\cup M_+^*\cup M_-\cup
M_-^*$
be a decomposition into disjoint subsets, let $\Bc$ be the distribution
\bqn
\Bc &=& b\prod_{r\in M_0} \Bom_{\alpha_r,\beta_r,\theta_r}
\prod_{r\in M_+\cup M_+^*}\Bet_{\gamma_r,\theta_r}
\prod_{r\in M_-\cup M_-^*}\Bxi_{\delta_r,\theta_r},
\eqn
and assume that the following conditions are satisfied:
\begin{itemize}
\item[(a)]
$b\in\G\Cti{}$;
\item[(b)]
$2\alpha_r\notin\Zm$, $\alpha_r+\beta_r\notin\Zmn$,
$\alpha_r-\beta_r\notin\Zmn$ for each $r\in M_0$;
\item[(c)]
$\gamma_r\notin\Zmn$ for each $r\in M_+$;
\item[(d)]
$\gamma_r\in\Zm$ for each $r\in M_+^*$;
\item[(e)]
$\delta_r\notin\Zmn$ for each $r\in M_-$;
\item[(f)]
$\delta_r\in\Zm$ for each $r\in M_-^*$.
\item[(g)]
$\ro_0<0$ {\rm (}or, equivalently, $\ro_1+\ro_2<0${\rm )}, where
\bqn
\ro_1 &=& \max\{-1-2\Re\beta_r\;:\;r\in M_0\}
\cup\{-1+\Re\delta_r\;:\;r\in M_-\},\\
\ro_2 &=& \max\{-1+2\Re\beta_r\;:\;r\in M_0\}
\cup\{-1+\Re\gamma_r\;:\;r\in M_+\},\\
\ro_0^* &=& \left\{\ba{ll}-1&\mbox{ if }M_0\neq\emptyset\\
-\iy&\mbox{ if }M_0=\emptyset,\ea\right.\\
\ro_0 &=& \max\{\ro_0^*,\ro_1+\ro_2\}.
\eqn
\end{itemize}
Finally, define the following distributions and constants:
\bqn
\Bc_+ &=& G[b]^{\frac{1}{2}}
b_+\prod_{r\in M_0\cup M_+\cup M_+^*}\Bet_{\gamma_r,\theta_r},
\quad
\Bc_- =   G[b]^{\frac{1}{2}}
b_-\prod_{r\in M_0\cup M_-\cup M_-^*}\Bxi_{\delta_r,\theta_r},\quad
\label{f5.16}
\\
\Omega_T &=& \sum_{r\in M_0} (\alpha_r^2-\beta_r^2),\label{f5.17}
\\
E_T &=& E[b]
\prod_{r\in M_0}\frac{G(1+\alpha_r+\beta_r)G(1+\alpha_r-\beta_r)}
{G(1+2\alpha_r)}
\prod_{\subthree{r\in M_0\cup M_+\cup M_+^*}%
{s\in M_0\cup M_-\cup M_-^*}{r\neq s}}
(1-t_st_r^{-1})^{-\gamma_r\delta_s}
\nn\\
&&\times
\prod_{r\in M_0\cup M_-\cup M_-^*} b_+(t_r)^{-\delta_r}
\prod_{r\in M_0\cup M_+\cup M_+^*} b_-(t_r)^{-\gamma_r},
\label{f5.18}
\eqn
where $\gamma_r=\alpha_r+\beta_r$ and $\delta_r=\alpha_r-\beta_r$
for $r\in M_0$.
Then
\begin{itemize}
\item[(i)]
$T_N(\Bc)$ is invertible for all sufficiently large $N$;
\item[(ii)]
$\Bc$ effects $\cR$-convergence with respect to
$[\Bc_+,\Bc_-]$ and $(\ros)$;
\item[(iii)]
$\det T_N(\Bc)=G[b]^NN^{\Omega_T} E_T(1+O(N^{\ro_0}))$.
\end{itemize}
\end{theorem}

Since $\Bxi_1\iv=\Bxi_{-1,0}$ and $\Bxi_{-1}=\Bxi_{1,\pi}$,
where $\Bxi_{\pm1}$ are the distributions defined in the paragraph before
Proposition \ref{p2.5}, the previous theorem implies that
the distribution $\Bh=\Bxi_1\iv\Bxi_{-1}$ effects $\cR$-convergence
with respect to $[1,\Bh]$ and $(-\iy,0,-\iy)$.
This is what we have used in the proof of Theorem \ref{t2.2}.

Now we specialize the above theorem to the case of
distributions which are even and which have no singularities at $-1$ and
$1$.
Note that
\bqn
\tilde{\Bom}_{\alpha_r,\beta_r,\theta_r}=\Bom_{\alpha_r,-\beta_r,-\theta_r}
\qquad
\tilde{\Bet}_{\gamma_r,\theta_r}=\Bxi_{\gamma_r,-\theta_r}.
\eqn
With this observation, it is easy to single out the class of distributions
of Fisher-Hartwig type which are even.

\begin{corollary}\label{c4.4}
Let $\theta_1,\dots,\theta_R\in(-\pi,0)\cup(0,\pi)$ be such that
$|\theta_r|\neq|\theta_s|$ for all $1\le r<s\le R$, $R\ge0$, put
$t_r=e^{i\theta_r}$, let $\{1,\dots,R\}=M_0\cup M_\pm\cup M_\pm^{*}$
be a decomposition into three disjoint subsets, let $\Bc$ be the
distribution
\bqn\label{f5.20}
\Bc &=& b\prod_{r\in M_0}
\Bom_{\alpha_r,\beta_r,\theta_r}
\Bom_{\alpha_r,-\beta_r,-\theta_r}
\prod_{r\in M_\pm\cup M_\pm^*}
\Bet_{\gamma_r,\theta_r}\Bxi_{\gamma_r,-\theta_r},
\eqn
and assume that the following conditions are satisfied:
\begin{itemize}
\item[(a)]
$b\in \G\Cti{}$ is even;
\item[(b)]
$2\alpha_r\notin\Zm$, $\alpha_r+\beta_r\notin\Zmn$,
$\alpha_r-\beta_r\notin\Zmn$ and $|\Re\beta_r|<1/2$
for each $r\in M_0$;
\item[(c)]
$\gamma_r\notin\Zmn$ and $\Re\gamma_r<1$ for each $r\in M_\pm$;
\item[(d)]
$\gamma_r\in\Zm$ for each $r\in M_\pm^*$.
\end{itemize}
Define the numbers
\bqn
\ro_{12} &=& \max\{-1+2|\Re\beta_r|\;:\;r\in M_0\}\cup
\{-1+\Re\gamma_r\;:\;r\in M_\pm\},\label{f5.21}\\
\ro_0^* &=& \left\{\ba{ll}-1&\mbox{ if }M_0\neq\emptyset\\
-\iy&\mbox{ if }M_0=\emptyset,\ea\right.\label{f5.22}\\
\ro_0 &=& \max\{\ro_0^*,2\ro_{12}\}, \label{f5.23}
\eqn
and the following distribution and constants:
\bqn
\Bc_+ &=&
G[b]^{\frac{1}{2}}b_+\prod_{r\in M_0}\Bet_{\alpha_r+\beta_r,\theta_r}
\Bet_{\alpha_r-\beta_r,-\theta_r}
\prod_{r\in M_\pm\cup M_\pm^*}
\Bet_{\gamma_r,\theta_r}\label{f4.25}
\\
\Omega_T^\sym &=& 2\sum_{r\in M_0} (\alpha_r^2-\beta_r^2),
\label{f4.26}
\eqn
\bqn
E_T^\sym &=& E[b]
\prod_{r\in M_0}\frac{G^2(1+\alpha_r+\beta_r)G^2(1+\alpha_r-\beta_r)}
{G^2(1+2\alpha_r)}
\nn\\
&&\times
\prod_{r\in M_0}
b_+(t_r)^{-2(\alpha_r-\beta_r)}b_-(t_r)^{-2(\alpha_r+\beta_r)}
\prod_{r\in M_\pm\cup M_\pm^*}
b_-(t_r)^{-2\gamma_r}
\nn\\
&&\times
\prod_{r,s\in M_0}
(1-t_rt_s)^{-(\alpha_r-\beta_r)(\alpha_s-\beta_s)}
(1-t_r\iv t_s\iv)^{-(\alpha_r+\beta_r)(\alpha_s+\beta_s)}
\nn\\
&&\times
\prod_{\subtwo{r,s\in M_0}{r\neq s}}
(1-t_rt_s\iv)^{-(\alpha_r-\beta_r)(\alpha_s+\beta_s)}
(1-t_r\iv t_s)^{-(\alpha_r+\beta_r)(\alpha_s-\beta_s)}
\nn\\
&&\times
\prod_{\subtwo{r\in M_0}{s\in M_\pm\cup M_\pm^*}}
(1-t_r\iv t_s\iv)^{-2(\alpha_r+\beta_r)\gamma_s}
(1-t_rt_s\iv)^{-2(\alpha_r-\beta_r)\gamma_s}
\nn\\
&&\times
\prod_{r,s\in M_\pm\cup M_\pm^*}
(1-t_r\iv t_s\iv)^{-\gamma_r\gamma_s}.
\label{f4.27}
\eqn
Then
\begin{itemize}
\item[(i)]
$T_N(\Bc)$ is invertible for all sufficiently large $N$;
\item[(ii)]
$\Bc$ effects $\cR$-convergence with respect to
$[\Bc_+,\tilde{\Bc}_+]$ and $(\ro_0,\ro_{12},\ro_{12})$;
\item[(iii)]
$\det T_N(\Bc)=G[b]^NN^{\Omega_T^\sym} E_T^\sym(1+O(N^{\ro_0}))$.
\end{itemize}
\end{corollary}
\begin{proof}
This corollary is just a special case of Theorem \ref{t4.3new}. The setting
in which
we have to apply Theorem \ref{t4.3new} is the following. The number $R$ of
Theorem
\ref{t4.3new} is twice the number $R$ of this corollary. More precisely,
the set $M_0$ of Theorem \ref{t4.3new} has to be identified with two copies
of the present set $M_0$. Let us denote these two copies by $M_0^{(1)}$ and
$M_0^{(2)}$. The set $M_+$ ($M_+^*$, resp.) of Theorem \ref{t4.3new} has to
be identified with one copy of the present set $M_\pm$ ($M_\pm^*$, resp.).
In the same way, the set $M_-$ ($M_-^*$, resp.) of Theorem \ref{t4.3new} has
also to
be identified with one copy of the present set $M_\pm$ ($M_\pm^*$, resp.).
The parameters corresponding to the index set
\begin{itemize}
\item
$M_0^{(1)}$ are $(\alpha_r,\beta_r,\theta_r)$ for $r\in M_0$;
\item
$M_0^{(2)}$ are $(\alpha_r,-\beta_r,-\theta_r)$ for $r\in M_0$;
\item
$M_+\cup M_+^*$ are $(\gamma_r,\theta_r)$ for $r\in M_\pm\cup M_\pm^*$;
\item
$M_-\cup M_-^*$ are $(\gamma_r,-\theta_r)$ for $r\in M_\pm\cup M_\pm^*$.
\end{itemize}
{}From this it is readily seen that $\ro_1=\ro_2=\ro_{12}$. In order for
$\ro_0$,
or equivalently $\ro_1+\ro_2=2\ro_{12}$, to be negative, it is necessary and
sufficient that $|\Re\beta_r|<1/2$ and $\Re\gamma_r<1$.
Hence conditions (a)--(d) of this corollary imply conditions (a)--(g) of
Theorem \ref{t4.3new}.

The distributions $\Bc_+$ and $\Bc_-$ defined in (\ref{f5.16})
are given by (\ref{f4.25}) and
\bqn
\Bc_- &=&
G[b]^{\frac{1}{2}}b_-\prod_{r\in M_0}\Bxi_{\alpha_r+\beta_r,-\theta_r}
\Bxi_{\alpha_r-\beta_r,\theta_r}
\prod_{r\in M_\pm\cup M_\pm^*}
\Bxi_{\gamma_r,-\theta_r}.
\eqn
Since $b$ is odd, we have $b_-=\tilde{b}_+$ and thus $\Bc_-=\tilde{\Bc}_+$.
Because of the above parameters corresponding to $M_0^{(1)}$ and
$M_0^{(2)}$,
it is easily seen that $\Omega_T$ as given in
(\ref{f5.17}) becomes $\Omega_T^\sym$ as given above.
It is somewhat troublesome, but straightforward to verify that
$E_T$ as given in (\ref{f5.18}) becomes $E_T^\sym$.
This completes the proof.
\end{proof}


\section{Asymptotics of determinants of symmetric \\ Toeplitz plus Hankel
matrices with \\ Fisher-Hartwig distributions}
\label{sec:6}

In this section we finally combine the Limit Theorem (Theorem \ref{t2.2})
and Theorem \ref{t2.3} with Corollary \ref{c4.4} in order to obtain an
asymptotic formula
for the determinants of symmetric Toeplitz + Hankel matrices.
This result is based on the asymptotic formula for
determinants of (symmetric) Toeplitz matrices.

In what follows, let $G[b]$ and $\wh{E}[b]$ be the constants (\ref{f1.Gb})
and (\ref{f1.Ehb}), and let $b_\pm$ be the functions (\ref{f1.bpm}).

\begin{theorem}\label{t5.1}
Let $\Bc$ be a distribution that fulfills the assumptions of
Corollary {\sl \ref{c4.4}}. Define the number $\ro_{12}$ by \refeq{f5.21}
and the constants
\bqn
\Omega_M^\sym &=&
\sum_{r\in M_0}\left(\alpha_r^2-\beta_r^2\right),
\\
E_M^\sym &=& \wh{E}[b]
\prod_{r\in M_0}\frac{G(1+\alpha_r+\beta_r)G(1+\alpha_r-\beta_r)}
{G(1+2\alpha_r)}
\nn\\
&&\times
\prod_{r\in M_0}
\frac{(1-t_r\iv)^{(\alpha_r+\beta_r)/2}(1-t_r)^{(\alpha_r-\beta_r)/2}}
{(1+t_r\iv)^{(\alpha_r+\beta_r)/2}(1+t_r)^{(\alpha_r-\beta_r)/2}}
2^{\alpha_r^2-\beta_r^2}
\!\!\!\!\!\!\prod_{r\in M_\pm\cup M_\pm^*}
\frac{(1-t_r\iv)^{\gamma_r/2}}{(1+t_r\iv)^{\gamma_r/2}}
\nn\\
&&\times
\prod_{r\in M_0}
b_+(t_r)^{-(\alpha_r-\beta_r)}
b_-(t_r)^{-(\alpha_r+\beta_r)}
\prod_{r\in M_\pm\cup M_\pm^*}
b_-(t_r)^{-\gamma_r}
\nn\\
&&\times
\prod_{r,s\in M_0}
(1-t_rt_s)^{-(\alpha_r-\beta_r)(\alpha_s-\beta_s)/2}
(1-t_r\iv t_s\iv)^{-(\alpha_r+\beta_r)(\alpha_s+\beta_s)/2}
\nn\\
&&\times
\prod_{\subtwo{r,s\in M_0}{r\neq s}}
(1-t_rt_s\iv)^{-(\alpha_r-\beta_r)(\alpha_s+\beta_s)/2}
(1-t_r\iv t_s)^{-(\alpha_r+\beta_r)(\alpha_s-\beta_s)/2}
\nn\\
&&\times
\prod_{\subtwo{r\in M_0}{s\in M_\pm\cup M_\pm^*}}
(1-t_r\iv t_s\iv)^{-(\alpha_r+\beta_r)\gamma_s}
(1-t_rt_s\iv)^{-(\alpha_r-\beta_r)\gamma_s}
\nn\\
&&\times
\prod_{r,s\in M_\pm\cup M_\pm^*}
(1-t_r\iv t_s\iv)^{-\gamma_r\gamma_s/2}.
\eqn
Then
\bqn
\det M_N(\Bc) &=&
G[b]^NN^{\Omega_M^\sym} E_M^\sym(1+O(N^{\ro_{12}})).
\eqn
\end{theorem}
\begin{proof}
The distribution $\Bc$ belongs to $\cD'(K\cup\wt{K})$,
where $K=\{t_r: 1\le r\le R\}$. Thus the assumptions of
Theorem \ref{t2.3} are fulfilled and we obtain
\bqn
\det T_{2N}(\Bsi \Bc) &=& (\det M_N(\Bc))^2.\nn
\eqn
Moreover, the distribution $\Bc$ has the smooth part
\bqn
c(e^{i\theta}) &=& b(e^{i\theta})
\!\prod_{r\in M_0}
\omega_{\alpha_r,\beta_r,\theta_r}(e^{i\theta})
\omega_{\alpha_r,-\beta_r,-\theta_r}(e^{i\theta})
\!\!\!\!\prod_{r\in M_\pm\cup M_\pm^*}
\eta_{\gamma_r,\theta_r}(e^{i\theta})
\xi_{\gamma_r,-\theta_r}(e^{i\theta}),\nn
\eqn
{}From Corollary \ref{c4.4} (ii) it follows that $\Bc$ effects
$\cR$-convergence with respect to $[\Bc_+,\tilde{\Bc}_+]$ and
$(\ro_0,\ro_{12},\ro_{12})$, where $\ro_0$ is defined by (\ref{f5.22}) and
(\ref{f5.23}) and $\Bc_+$ is the distribution given by (\ref{f4.25}).
Notice that $\ro_0$ and $\ro_{12}$ are negative real numbers.
Obviously, $\Bc_+$ belongs to
$\cG\Dp(K\cup\wt{K})$ and has the smooth part
\bqn
c_+(e^{i\theta}) &=&
G[b]^{\frac{1}{2}}b_+(e^{i\theta})\prod_{r\in M_0}
\eta_{\alpha_r+\beta_r,\theta_r}(e^{i\theta})
\eta_{\alpha_r-\beta_r,-\theta_r}(e^{i\theta})
\prod_{r\in M_\pm\cup M_\pm^*}
\eta_{\gamma_r,\theta_r}(e^{i\theta}).\nn
\eqn
It is readily seen that $c(t)=c_+(t)\tilde{c}_+(t)$. Hence the assumptions
of Theorem \ref{t2.2} are fulfilled with $\Ba$ and $\Ba_+$ replaced by
$\Bc$ and $\Bc_+$, respectively, and $K$ replaced by $K\cup\wt{K}$.
Thus
\bqn
\det T_{2N}(\Bsi \Bc) &=& \det T_{2N}(\Bc)\left(
\frac{c_+(1)}{c_+(-1)}+O\left(N^{\max\{\ro_0,\ro_{12}\}}\right)\right).\nn
\eqn
Note that $\max\{\ro_0,\ro_{12}\}=\max\{\ro_0^*,2\ro_{12},\ro_{12}\}=
\max\{\ro_0^*,\ro_{12}\}=\ro_{12}$.
Combining the previous formulas we obtain
\bqn
(\det M_N(\Bc))^2 &=&
\det T_{2N}(\Bc)\frac{c_+(1)}{c_+(-1)}
\left(1+O(N^{\ro_{12}})\right).\nn
\eqn
The asymptotics of $\det T_{2N}(\Bc)$ follows from
Corollary \ref{c4.4} (iii). Observe the change from $N$ to $2N$ and that
$\ro_0\le\ro_{12}$. Thus
\bqn
(\det M_N(\Bc))^2 &=&
G[b]^{2N}(2N)^{\Omega^\sym_T}E^\sym_T\frac{c_+(1)}{c_+(-1)}
\left(1+O\left(N^{\ro_{12}}\right)\right),\nn
\eqn
where $\Omega^\sym_T$ and $E^\sym_T$ are the constants
(\ref{f4.26}) and (\ref{f4.27}). Notice also that
\bqn
2^{\Omega^\sym_T}\frac{c_+(1)}{c_+(-1)}
&=&
\prod_{r\in M_0}
2^{2(\alpha_r^2-\beta_r^2)}
\frac{(1-t_r\iv)^{\alpha_r+\beta_r}(1-t_r)^{\alpha_r-\beta_r}}%
{(1+t_r\iv)^{\alpha_r+\beta_r}(1+t_r)^{\alpha_r-\beta_r}}
\nn\\&&\times
\frac{b_+(1)}{b_+(-1)}
\prod_{r\in M_\pm\cup M_\pm^*}
\frac{(1-t_r\iv)^{\gamma_r}}{(1+t_r\iv)^{\gamma_r}}.\nn
\eqn
Taking into account that
\bqn
\wh{E}[b]^2 &=& E[b]\frac{b_+(1)}{b_+(-1)},\nn
\eqn
it is readily seen that
\bqn
2^{\Omega^\sym_T}E^\sym_T\frac{c_+(1)}{c_+(-1)}
&=& (E^\sym_M)^2.\nn
\eqn
Obviously, $\Omega^\sym_T=2\Omega^\sym_M$. Hence the last asymptotic
formula becomes
\bqn\label{f6.4}
(\det M_N(\Bc))^2 &=&
G[b]^{2N}N^{2\Omega^\sym_M}(E^\sym_M)^2\left(1+O(N^{\ro_{12}})\right).
\eqn
{}From this the desired asymptotic formula follows, up to a sign,
which will be determined by the following argument.
Let $U\subset \C^{d}$, where $d=2|M_0|+|M_\pm\cup M_\pm^*|$,
be the set of all $d$-tuples
\begin{equation}\label{f.param}
z=\left((\alpha_r,\beta_r)_{r\in M_0},
(\gamma_r)_{r\in M_\pm}, (\gamma_r)_{r\in M_\pm^*}\right)
\end{equation}
such that
\begin{itemize}
\item[(a)]
$2\alpha_r\notin\Zm$, $\alpha_r\pm\beta_r\notin\Zm\cup\{0\}$,
$|\Re\beta_r|<1/2$ for all $r\in M_0$;
\item[(b)]
$\gamma_r\notin\Zm\cup\{0\}$, $\Re\gamma_r<1$
for all $r\in M_\pm\cup M_\pm^*$.
\end{itemize}
Let the even function $b_1\in G\Cti{}$ be arbitrary but fixed, and introduce
$$
f_N(z)=\frac{\det M_N(\Bc)}{G[b]^NN^{\Omega^\sym_M}},\quad
f(z)=E^\sym_M,
$$
where $\Bc$ is the distribution (\ref{f5.20}) with the parameters
given by (\ref{f.param}), and the constants $\Omega^\sym_M$ and
$E^\sym_M$ are defined correspondingly.

{}From what has been proved so far, it follows that
\bqn
(f_N(z))^2\to (f(z))^2
\eqn
for $z\in U$. Moreover, the convergence is uniform on compact
subsets of $U$. Since $U$ is connected and since $f(z)\neq0$ for all
$z\in U$, it follows that either $f_N(z)\to f(z)$ on all of $U$ or
$f_N(z)\to -f(z)$ on all of $U$.

Now let $U_0$ stand for the set of all parameters (\ref{f.param})
for which
\begin{itemize}
\item[(a)]
$2\alpha_r\notin\Zm$, $|\Re\beta_r|<1/2$ for all $r\in M_0$;
\item[(b)]
$\Re\gamma_r<1$
for all $r\in M_\pm\cup M_\pm^*$.
\end{itemize}
Because $f_N(z)$ and $f(z)$ are functions that depend analytically
on $z$, because of the uniform convergence on compact subsets of $U$
and because of the concrete structure of $U$ and $U_0$, it follows
that either $f_N(z)\to f(z)$ on all of $U_0$ or
$f_N(z)\to -f(z)$ on all of $U_0$.

For $z=0$ we know that $0\in U_0$ and that
$f_N(0)\to f(0)$. Thus we can conclude that
$f_N(z)\to f(z)$ on all of $U_0$.
Hence the desired asymptotic formula with the correct sign follows.
\end{proof}


\section{Asymptotics of determinants of symmetric \\ Toeplitz plus Hankel
matrices with piecewise continuous functions}
\label{sec:7}

{}From Theorem \ref{t5.1} we obtain immediately the following result
concerning the asymptotics of the determinants of symmetric
Toeplitz plus Hankel matrices $M_N(\phi)$ with a particular piecewise
continuous
generating function.

\begin{theorem}\label{t7.1}
Let $\theta_0\in(0,\pi)$ and $\beta\in\C$ be such that $|\Re\beta|<1/2$.
Put $t_0=e^{i\theta_0}$. Then
\bqn
\lim_{N\to\iy}
\frac{\det M_N(t_{\beta,\theta_0}t_{-\beta,-\theta_0})}
{N^{-\beta^2}} &=& E,\nn
\eqn
where
\bqn
E &=&
2^{-\beta^2}
(1-t_0^2)^{-\frac{\beta^2}{2}}(1-t_0^{-2})^{-\frac{\beta^2}{2}}
\frac{(1-t_0\iv)^{\frac{\beta}{2}}(1-t_0)^{-\frac{\beta}{2}}}%
{(1+t_0\iv)^{\frac{\beta}{2}}(1+t_0)^{-\frac{\beta}{2}}}
G(1+\beta)G(1-\beta) .\nn
\eqn
\end{theorem}
\begin{proof}
In Theorem \ref{t5.1} we put $M_\pm =M_\pm^*=\emptyset$, $M_0=\{1\}$,
$\alpha_1=0$, $\beta_1=\beta$, $\theta_1=\theta_0$ and $b(t)=1$. Observe
that
$\Bom_{0,\beta,\theta_0}\Bom_{0,-\beta,-\theta_0}=
t_{\beta,\theta_0}t_{-\beta,-\theta_0}$.
\end{proof}

The main result concerning the asymptotics of the determinants of symmetric
Toeplitz + Hankel matrices $M_N(\phi)$ with ``general'' piecewise continuous
generating functions is the following theorem.
\begin{theorem}\label{t7.2}
Let
\bqn
c(e^{i\theta}) &=&
b(e^{i\theta})\prod_{r=1}^R
t_{\beta_r,\theta_r}(e^{i\theta})t_{-\beta_r,-\theta_r}(e^{i\theta})
\eqn
where $b\in \G\B$ is an even function, $\theta_1,\dots,\theta_R\in(0,\pi)$
are distinct numbers, and $\beta_1,\dots,\beta_R\in\C$ are such that
$|\Re\beta_r|<1/2$ for all $1\le r\le R$.
Let $G[b]$ and $\wh{E}[b]$ be the constants
\refeq{f1.Gb} and \refeq{f1.Ehb}, $b_\pm$ be the functions \refeq{f1.bpm},
$t_r=e^{i\theta_r}$, $1\le r\le R$, and introduce the constants
\bqn
\Omega_M^\sym &=&
-\sum_{r=1}^R\beta_r^2,
\\
E_M^\sym &=&
\wh{E}[b]
\prod_{r=1}^R
G(1+\beta_r)G(1-\beta_r)
(1-t_r^2)^{-\beta_r^2/2}
(1-t_r^{-2})^{-\beta_r^2/2}
\nn\\
&&\times
\!\!\!\prod_{1\le r<s\le R}\!\!\!\!\!
(1-t_rt_s)^{-\beta_r\beta_s}
(1-t_r\iv t_s\iv)^{-\beta_r\beta_s}
(1-t_rt_s\iv)^{\beta_r\beta_s}
(1-t_r\iv t_s)^{\beta_r\beta_s}
\nn\\
&&\times
\prod_{r=1}^R
2^{-\beta_r^2}
\frac{(1-t_r\iv)^{\beta_r/2}(1-t_r)^{-\beta_r/2}}
{(1+t_r\iv)^{\beta_r/2}(1+t_r)^{-\beta_r/2}}
\prod_{r=1}^R
b_+(t_r)^{\beta_r}b_-(t_r)^{-\beta_r}.
\eqn
 Then
\bqn
\lim_{N\to\iy}
\frac{\det M_N(\Bc)}{G[b]^NN^{\Omega_M^\sym}} &=& E_M^\sym.
\eqn
\end{theorem}
\begin{proof}
The asymptotic formula follows from Corollary \ref{c2.5new}),
Theorem \ref{t7.1} and the asymptotic formula (\ref{f2.7new}).
\end{proof}

We remark that in the case $b\in G_1 C^\iy(\T)$, the previous theorem
follows also directly from Theorem \ref{t5.1}.

In view of the asymptotic formulas established in
Theorem \ref{t2.3new}, Theorem \ref{t2.4new} and Theorem \ref{t7.1}
we now raise the following conjecture.

\begin{conjecture}
Let $\theta_0\in(0,\pi)$ and $\beta_1,\beta_2\in\C$ be such that
$|\Re\beta_1|<1/2$ and $|\Re\beta_2|<1/2$ and $|\Re\beta_1+\Re\beta_2|<1/2$.
Put $t_0=e^{i\theta_0}$.
Then
\bqn
\lim_{N\to\iy}
\frac{\det M_N(t_{\beta_1,\theta_0}t_{\beta_2,-\theta_0})}{N^\Omega}=E,
\eqn
where
\bqn
\Omega &=& -\beta_1^2-\beta_1\beta_2-\beta_2^2,
\\
E &=& G(1+\beta_1)G(1+\beta_2)G(1-\beta_1-\beta_2)
2^{\beta_1\beta_2}
\nn\\
&&\times
(1-t_0^{-2})^{\beta_1^2/2+\beta_1\beta_2}
(1-t_0^{2})^{\beta_2^2/2+\beta_1\beta_2}
\nn\\
&&\times
\frac{(1-t_0\iv)^{\beta_1/2}(1-t_0)^{\beta_2/2}}
{(1+t_0\iv)^{\beta_1/2}(1+t_0)^{\beta_2/2}}.
\eqn
\end{conjecture}
This conjecture fits with the results established in Theorem \ref{t2.3new},
Theorem \ref{t2.4new} and Theorem \ref{t7.1}.


\end{document}